
\def\LaTeX{%
  \let\Begin\begin
  \let\End\end
  \let\salta\relax
  \let\finqui\relax
  \let\futuro\relax}

\def\UK{\def\our{our}\let\sz s}
\def\USA{\def\our{or}\let\sz z}

\UK



\LaTeX

\USA


\salta

\documentclass[twoside,12pt]{article}
\setlength{\textheight}{24cm}
\setlength{\textwidth}{16cm}
\setlength{\oddsidemargin}{2mm}
\setlength{\evensidemargin}{2mm}
\setlength{\topmargin}{-15mm}
\parskip2mm


\usepackage[usenames,dvipsnames]{color}
\usepackage{amsmath}
\usepackage{amsthm}
\usepackage{amssymb}
\usepackage[mathcal]{euscript}


\usepackage[T1]{fontenc}
\usepackage[latin1]{inputenc}
\usepackage[english]{babel}
\usepackage[babel]{csquotes}

\usepackage{cite}

\usepackage{latexsym}
\usepackage{graphicx}
\usepackage{mathrsfs}
\usepackage{mathrsfs}
\usepackage{hyperref}
\usepackage{pgfplots}

%
%


\definecolor{viola}{rgb}{0.3,0,0.7}
\definecolor{ciclamino}{rgb}{0.5,0,0.5}
\definecolor{bluette}{rgb}{0,0.4,0.8}

\def\newjuerg #1{{\color{red}#1}}
\def\pier #1{{\color{blue}#1}}

\def\newjuerg #1{#1}
\def\pier #1{#1}






\bibliographystyle{plain}


%

\finqui

\def\Beq{\Begin{equation}}
\def\Eeq{\End{equation}}
\def\Bsist{\Begin{eqnarray}}
\def\Esist{\End{eqnarray}}

\def\Bthm{\Begin{theorem}}
\def\Ethm{\End{theorem}}
\def\Blem{\Begin{lemma}}
\def\Elem{\End{lemma}}

\def\Brem{\Begin{remark}\rm}
\def\Erem{\End{remark}}

\def\Bdim{\Begin{proof}}
\def\Edim{\End{proof}}
\def\Bcenter{\Begin{center}}
\def\Ecenter{\End{center}}
\let\non\nonumber




\def\step #1 \par{\medskip\noindent{\bf #1.}\quad}


\def\Lip{Lip\-schitz}

\def\aand{\quad\hbox{and}\quad}

\def\lhs{left-hand side}
\def\rhs{right-hand side}


\def\bhv{behavi\our}


\def\multibold #1{\def\arg{#1}%
  \ifx\arg\pto \let\next\relax
  \else
  \def\next{\expandafter
    \def\csname #1#1#1\endcsname{{\bf #1}}%
    \multibold}%
  \fi \next}

\def\pto{.}

\def\multical #1{\def\arg{#1}%
  \ifx\arg\pto \let\next\relax
  \else
  \def\next{\expandafter
    \def\csname cal#1\endcsname{{\cal #1}}%
    \multical}%
  \fi \next}


\def\multimathop #1 {\def\arg{#1}%
  \ifx\arg\pto \let\next\relax
  \else
  \def\next{\expandafter
    \def\csname #1\endcsname{\mathop{\rm #1}\nolimits}%
    \multimathop}%
  \fi \next}

\multibold
qwertyuiopasdfghjklzxcvbnmQWERTYUIOPASDFGHJKLZXCVBNM.

\multical
QWERTYUIOPASDFGHJKLZXCVBNM.

\multimathop
diag dist div dom mean meas sign supp .


\def\accorpa #1#2{\eqref{#1}--\eqref{#2}}
\def\Accorpa #1#2 #3 {\gdef #1{\eqref{#2}--\eqref{#3}}%
  \wlog{}\wlog{\string #1 -> #2 - #3}\wlog{}}


\def\separa{\noalign{\allowbreak}}

\def\somma #1#2#3{\sum_{#1=#2}^{#3}}

\def\<#1>{\mathopen\langle #1\mathclose\rangle}
\def\norma #1{\mathopen \| #1\mathclose \|}

\def\[#1]{\mathopen\langle\!\langle #1\mathclose\rangle\!\rangle}

\def\iot {\int_0^t}
\def\ioT {\int_0^T}
\def\intQt{\int_{Q_t}}
\def\intQ{\int_Q}
\def\iO{\int_\Omega}

\def\dt{\partial_t}

\def\cpto{\,\cdot\,}

\def\checkmmode #1{\relax\ifmmode\hbox{#1}\else{#1}\fi}
\def\aeO{\checkmmode{a.e.\ in~$\Omega$}}
\def\aeQ{\checkmmode{a.e.\ in~$Q$}}

\def\aet{\checkmmode{a.e.\ in~$(0,T)$}}

\def\aat{\checkmmode{for a.a.~$t\in(0,T)$}}


\def\erre{{\mathbb{R}}}

\def\enne{{\mathbb{N}}}




\def\genspazio #1#2#3#4#5{#1^{#2}(#5,#4;#3)}
\def\spazio #1#2#3{\genspazio {#1}{#2}{#3}T0}

\def\L {\spazio L}
\def\H {\spazio H}

\def\C #1#2{C^{#1}([0,T];#2)}


\def\Lx #1{L^{#1}(\Omega)}
\def\Hx #1{H^{#1}(\Omega)}

\def\LQ #1{L^{#1}(Q)}

\def\Luno{\Lx 1}
\def\Ldue{\Lx 2}

\def\Huno{\Hx 1}
\def\Hdue{\Hx 2}


\def\LQ #1{L^{#1}(Q)}


\let\theta\vartheta
\let\eps\varepsilon
\let\bphi\phi
\let\phi\varphi

\let\hat\widehat

\let\TeXchi\chi                         
\newbox\chibox
\setbox0 \hbox{\mathsurround0pt $\TeXchi$}
\setbox\chibox \hbox{\raise\dp0 \box 0 }
\def\chi{\copy\chibox}


\def\QED{\hfill $\square$}


\def\VA #1{V_A^{#1}}
\def\VB #1{V_B^{#1}}
\def\Az #1{A_0^{#1}}
\def\Vz #1{V_0^{#1}}

\def\Beta{\hat\beta}

\def\Pi{\hat\pi}
\def\Lpi{L_\pi}

\def\phiz{\phi_0}
\def\mz{m_0}
\def\1{{\bf 1}}

\def\phis{\phi_{\!\sigma}}
\def\mus{\mu_\sigma}
\def\fs{f_{\!\sigma}}
\def\betal{\beta_\lambda}
\def\Betal{\Beta_\lambda}
\def\phisl{\phis^\lambda}
\def\musl{\mus^\lambda}
\def\Asj{A^\sigma_j(\psi)}
\def\Azj{A^0_j(\psi)}
\def\rhos{\rho_\sigma}
\def\Alpha{\hat\alpha}

\def\betaeps{\beta_\eps}
\def\Betaeps{\Beta_\eps}
\let\gam\gamma
\def\Gam{\hat\gam}
\def\Vp{V^*}
\def\lims{\lim_{\sigma\searrow0}}
\def\liminfs{\liminf_{\sigma\searrow0}}
\def\phil{\bphi_\lambda}
\def\bphiz{\bphi_0}
\def\vn{v_n}

\Begin{document}


\title{An asymptotic analysis
  \\ for a generalized Cahn--Hilliard system
  \\ with fractional operators}
\author{}
\date{}
\maketitle
\Bcenter
\vskip-1.5cm
{\large\sc Pierluigi Colli$^{(1)}$}\\
{\normalsize e-mail: {\tt pierluigi.colli@unipv.it}}\\[.25cm]
{\large\sc Gianni Gilardi$^{(1)}$}\\
{\normalsize e-mail: {\tt gianni.gilardi@unipv.it}}\\[.25cm]
{\large\sc J\"urgen Sprekels$^{(2)}$}\\
{\normalsize e-mail: {\tt sprekels@wias-berlin.de}}\\[.45cm]
$^{(1)}$
{\small Dipartimento di Matematica ``F. Casorati''}\\
{\small Universit\`a di Pavia}\\
{\small via Ferrata 5, 27100 Pavia, Italy}\\[.45cm]
$^{(2)}$
{\small Department of Mathematics}\\
{\small Humboldt-Universit\"at zu Berlin}\\
{\small Unter den Linden 6, 10099 Berlin, Germany}\\[1.5mm]
{\small and}\\[1.5mm]
{\small Weierstrass Institute for Applied Analysis and Stochastics}\\
{\small Mohrenstrasse 39, 10117 Berlin, Germany}
\Ecenter
\Begin{abstract}\noindent
In the recent paper `Well-posedness and regularity for a generalized fractional Cahn--Hilliard system'
(Atti Accad.\ Naz.\ Lincei Rend.\ Lincei Mat.\ Appl.\ 30 (2019), 437-478), the same authors have studied 
viscous and nonviscous Cahn--Hilliard systems of two operator equations in which nonlinearities of
double-well type, like regular or logarithmic potentials, as well as nonsmooth potentials with 
indicator functions, were admitted. 
The operators appearing in the system equations are fractional powers $A^{2r}$ and $B^{2\sigma}$ (in~the spectral sense) 
of~general linear operators $A$ and~$B$, which are densely defined, 
unbounded, selfadjoint, and monotone in the Hilbert space $L^2(\Omega)$, for some bounded 
and smooth domain $\Omega\subset\erre^3$, and have compact resolvents.
Existence, uniqueness, and regularity results have been proved in the quoted paper. 
Here, in the case of the viscous system, we analyze the asymptotic behavior of the solution as the parameter 
$\sigma$ appearing in the operator $B^{2\sigma}$ decreasingly tends to zero. 
We prove convergence to a phase relaxation problem at the limit,
and we also investigate this limiting problem,
in which an additional term containing the projection of the phase variable on the kernel of $B$ appears.
\vskip3mm
\noindent {\bf Key words:}
Fractional operators, Cahn--Hilliard systems, asymptotic analysis. 
\vskip3mm
\noindent {\bf AMS (MOS) Subject Classification:} 35K45, 35K90, 35R11, 35B40.
\End{abstract}
\salta
\pagestyle{myheadings}
\newcommand\testopari{\sc Colli \ --- \ Gilardi \ --- \ Sprekels}
\newcommand\testodispari{\sc Generalized fractional Cahn--Hilliard system}
\markboth{\testopari}{\testodispari}
\finqui
%

\section{Introduction}
\label{Intro}
\setcounter{equation}{0}

A research project that the three of us recently carried out in \cite{CGS18,CGS21bis,CGS19},
deals with the well-posedness, regularity and optimal control for the abstract evolutionary system
\Bsist
  && \dt\phi + A^{2r} \mu = 0,
  \label{Iprima}
  \\
  && \tau \dt\phi + B^{2\sigma}\phi + F'(\phi) = \mu + f,
  \label{Iseconda}
  \\
  && \phi(0) = \phiz,
  \label{Icauchy}
\Esist
where $A^{2r}$ and $B^{2\sigma}$, with $r>0$ and $\sigma>0$,
denote fractional powers of the linear operators $A$ and~$B$, respectively.
These operators are supposed to be densely defined in $H:=\Ldue$, with $\Omega\subset\erre^3$,
selfadjoint and monotone, and to have compact resolvents. 
The above system is a generalization
of the standard or viscous Cahn--Hilliard system
(depending on whether $\tau=0$ or $\tau>0$), 
which models a phase separation process taking place in the container~$\Omega$.
The particular sample case $A^{2r}=B^{2\sigma}=-\Delta$ with homogeneous 
Neumann boundary conditions is included, indeed.
The physical variables $\,\phi\,$ and $\,\mu\,$ stand for the \emph{order parameter}
and the \emph{chemical potential}, respectively,
while $f$ is a given source term. 
Moreover, $\,F\,$ denotes a double-well potential. We offer three physically significant examples for $F$, namely,
\begin{align}
  & F_{reg}(r)
  := \frac 14 \, (r^2-1)^2 \,,
  \quad r \in \erre, 
  \label{regpot}
  \\[2mm]
  & F_{log}(r)
  := \left\{\begin{array}{ll}
    (1+r)\ln (1+r)+(1-r)\ln (1-r) - c_1 r^2\,,
    & \quad r \in (-1,1)
    \\[1mm]
    2\,\newjuerg{\ln}(2)-c_1\,,
    & \quad r\in\{-1,1\}
    \\[1mm]
    +\infty\,,
    & \quad r\not\in [-1,1]
  \end{array}\right. ,
  \label{logpot}
  \\[2mm]
  & F_{2obs}(r) := - c_2 r^2 
  \quad \hbox{if $|r|\leq1$}
  \aand
  \newjuerg{F_{2obs}}(r) := +\infty
  \quad \hbox{if $|r|>1$},
  \label{obspot}
\end{align}

\noindent where the constants $c_i$ in \eqref{logpot} and \eqref{obspot} satisfy
$c_1>1$ and $c_2>0$, so that $F_{log}$ and $F_{2obs}$ are nonconvex.
These potentials are called the {\em classical regular potential}, the {\em logarithmic potential\/},
and the {\em double obstacle potential\/}, respectively.
In irregular situations like \eqref{obspot}, one has to split $F$ into a nondifferentiable convex part~$\Beta$ 
(the~indicator function of $[-1,1]$, in the case of \eqref{obspot}) and a smooth perturbation~$\Pi$.
At the same time, one has to replace the derivative of the convex part
by the subdifferential and to interpret \eqref{Iseconda} as a differential inclusion
or, equivalently, as a variational inequality involving $\Beta$ rather than its subdifferential, 
as actually done in~\cite{CGS18}.

Fractional versions of the Cahn--Hilliard system have been considered by different authors
and are the subject of several papers. 
As for references regarding well-posedness and related problems, 
a rather large list of citations is given in~\cite{CGS18}; 
we recall some concerned and recent literature also here, by mentioning \cite{AM, AkSS1, BMH, Gal2, YLZ, WCW}. 
Moreover, one can find a number of results regarding the asymptotic behavior of solutions,
for the standard Cahn--Hilliard equations, for variants \newjuerg{thereof}, and for systems 
including the Cahn--Hilliard equations: without any claim of completeness,  we can quote, e.g.,
\cite{AkSS2, BCST, CRW, ChFaPr, ClGu, CF, CGS17, CGS24, CGS22, Gal3, GS, GiMiSchi2, GrPeSch, WaWu, WuZh, ZhaoLiu, Zheng}.
These works mainly deal with the asymptotics with respect to parameters, or the study 
of the trajectories and related topics, or the existence of global or exponential 
attractors and their properties. A special role in our citations is played by the paper~\cite{CG1},
where the longtime \bhv\ of the solutions as well as an asymptotic analysis similar to the one we 
address here are investigated for a fractional system involving the Allen--Cahn equation.

In this paper, we consider the viscous case $\tau>0$ within the system \accorpa{Iprima}{Icauchy} and 
study the asymptotic \bhv\ of the solution as the parameter $\sigma$ involved in the operator $B^{2\sigma}$ tends to zero.
In this analysis, a crucial role is played by the orthogonal projection operator $P:H\to H$ on the kernel $\ker B$ of~$B$.
Indeed, if $(\phis,\mus)$ denotes the solution to system \accorpa{Iprima}{Icauchy} for an arbitrary $\sigma>0$,
we prove that $(\phis,\mus)$ converges
as $\sigma\searrow0$ to a solution $(\phi,\mu)$ to the system
\Bsist
  && \dt\phi + A^{2r} \mu = 0,
  \label{Iprimaz}
  \\
  && \tau \dt\phi + \phi - P\phi + F'(\phi) = \mu + f,
  \label{Isecondaz}
  \\
  && \phi(0) = \phiz.
  \label{Icauchyz}
\Esist
In general, the convergence occurs along a subsequence, but in the case when the 
limit pair $(\phi,\mu)$ uniquely solves \accorpa{Iprimaz}{Icauchyz}, then \newjuerg{the whole 
family $(\phis,\mus)$} converges to $(\phi,\mu)$ in the sense made precise by the 
statement of Theorem~\ref{Asymptotics} below. 
Moreover, let us point out that the 
component $\phi$ of the pair $(\phi,\mu)$ is always uniquely determined, as it follows 
from the continuous dependence result given by Theorem~\ref{Partuniqz}. 
In the last part of the paper, we also discuss the limiting problem \pier{by proving
a class of regularity results, quite interesting in our opinion, for which we have
to use some sophisticated tools of interpolation theory. 
Our approach may be considered as an extension and further investigation with respect to the asymptotic results of \cite[Section~7]{CG1},} in which a phase relaxation problem is obtained at the limit. 
Also in the present paper the equation~\eqref{Isecondaz} can be seen as an ordinary 
differential equation, but with a nonlocal structure due to the presence of the 
projection operator~$P$. 
Our contribution \pier{here} gives account of a new line of investigation that in our \newjuerg{opinion should} be further explored.

The rest of the paper is organized as follows:
in the next Section~\ref{STATEMENT}, we list our assumptions and state our results.
The corresponding proofs are given in the last two Sections~\ref{ASYMPTOTICS} and~\ref{LIMITING}.


\section{Statement of the problem and results}
\label{STATEMENT}
\setcounter{equation}{0}

In this section, we state precise assumptions and notations and present our results.
Our framework is the same as in~\cite{CGS18},
and we briefly recall it here, for the reader's convenience.
First of all, the open set $\Omega\subset\erre^3$ is assumed to be bounded, connected and smooth.
We use the notation
\Beq
  H := \Ldue
  \label{defH}
\Eeq
and denote by $\norma\cpto$ and $(\cpto,\cpto)$ the standard norm and inner product of~$H$.
As for the operators involved in our system, we postulate that
\Bsist
  && A:D(A)\subset H\to H
  \aand
  B:D(B)\subset H\to H
  \quad \hbox{are}
  \non
  \\
  && \hbox{unbounded, monotone, selfadjoint, linear operators with compact resolvents.} 
  \qquad
  \label{hpAB} 
\Esist

We denote by 
$\{\lambda_j\}$ and $\{\lambda'_j\}$ the nondecreasing sequences of the eigenvalues of $\,A\,$ and $\,B$,
and by $\{e_j\}$ and $\{e'_j\}$ the (complete) systems of the corresponding orthonormal eigenvectors,
that~is,
\begin{align}
  & A e_j = \lambda_j e_j, \quad
  B e'_j = \lambda'_j e'_j,
  \aand
  (e_i,e_j) = (e'_i,e'_j) = \delta_{ij}
  \quad \hbox{for $i,j=1,2,\dots$},
  \label{eigen}
  \\
  & 0 \leq \lambda_1 \leq \lambda_2 \leq \dots
  \aand
  0 \leq \lambda'_1 \leq \lambda'_2 \leq \dots,
  \label{eigenvalues}
  \hbox{ with } \,
  \lim_{j\to\infty} \lambda_j
  = \lim_{j\to\infty} \lambda'_j
  = + \infty {,}
\end{align}
where $\delta_{ij}$ denotes the Kronecker index. 
The power $A^r$ of $A$ with an arbitrary positive real exponent~$r$ is given~by
\Bsist
  && A^r v = \somma j1\infty \lambda_j^r (v,e_j) e_j
  \quad \hbox{for $v\in\VA r$},
  \quad \hbox{where}
  \label{defAr}
  \\
  && \VA r := D(A^r)
  = \Bigl\{ v\in H:\ \somma j1\infty |\lambda_j^r (v,e_j)|^2 < +\infty \Bigr\}.
  \label{defdomAr}
\Esist
In principle, we could endow $\VA r$ with the standard graph norm
in order to make $\VA r$ a Hilbert space.
However, we will choose an equivalent Hilbert structure later on.
In the same way, for $\sigma>0$,
we define the power $B^\sigma$ of~$B$.
For its domain we use the notation
\Bsist
  && \VB\sigma := D(B^\sigma),
  \quad \hbox{with the norm $\norma\cpto_{B,\sigma}$ associated to the inner product}
  \label{defBs}
  \non
  \\
  && (v,w)_{B,\sigma} := (v,w) + (B^\sigma v,B^\sigma w)
  \quad \hbox{for $v,w\in \VB\sigma$}.
  \label{defprodBs}
\Esist

At this point, we can start listing our assumptions.
First of all,
\Beq
  \hbox{$r$, $\sigma_0$ and $\tau$ are fixed positive numbers, and $\sigma\in(0,\sigma_0)$ is a parameter.}
  \label{hprst}
\Eeq
As for the linear operators, we postulate, besides~\eqref{hpAB}, that
\Bsist
  && \hbox{either} \quad
  \lambda_1 > 0 
  \quad \hbox{or} \quad
  \hbox{$0=\lambda_1<\lambda_2$ and $e_1$ is a constant;}
  \label{hpsimple}
  \\
  && \hbox{if\quad $\lambda_1=0$,\quad
  \newjuerg{then} the constant functions belong to $\VB\sigma$}.
  \label{hpVB}
\Esist
In \cite{CGS18} some remarks are given on the above assumptions.
Moreover, it is shown that an equivalent Hilbert structure on $\VA r$ is obtained
by taking the norm defined~by
\Beq
  \norma v_{A,r}^2 := \left\{ 
  \begin{aligned}
  & \norma{A^r v}^2
  = \somma j1\infty |\lambda_j^r (v,e_j)|^2
  \qquad \hbox{if $\lambda_1>0$,}
  \\
  & |(v,e_1)|^2 + \norma{A^r v}^2
  = |(v,e_1)|^2 + \somma j2\infty |\lambda_j^r (v,e_j)|^2
  \qquad \hbox{if $\lambda_1=0$,}
  \end{aligned}
  \right.
  \label{defnormaAr}
\Eeq
and the corresponding inner product, which we term $(\cpto,\cpto)_{A,r}$.
This equivalence is trivial if $\lambda_1>0$.
In the opposite case $\lambda_1=0$, with the notation
\Beq
  \mean v := \frac 1 {|\Omega|} \iO v 
  \qquad \hbox{for $v\in\Luno$}
  \label{defmean}
\Eeq
for the mean value ot the generic function~$v$,
the equivalence relies on the inequality
\Beq
  \norma v \leq C_P \, \norma{A^r v}
  \quad \hbox{for every $v\in\VA r$ with $\mean v=0$ \ if $\lambda_1=0$},
  \label{poincare}
\Eeq
which is of Poincar\'e type, 
since the term $(v,e_1)$ appearing in \eqref{defnormaAr} and involving the constant function~$e_1$ 
(see \eqref{hpsimple}) is proportional to~$\mean v$.
Next, the nonlinear potential $F$ appearing in \eqref{Iseconda} is split as follows:
\Bsist
  && F = \Beta + \Pi,
  \quad \hbox{where}
  \label{hpF}
  \\
  && \Beta : \erre \to [0,+\infty]
  \quad \hbox{is convex, proper and l.s.c.\ with} \quad
  \Beta(0) = 0;
  \label{hpBeta}
  \\
  \separa
  && \Pi : \erre \to \erre
  \quad \hbox{is of class $C^1$ with a \Lip\ continuous first derivative}; \qquad
  \label{hpPi}
  \\
  && \hbox{it holds} \quad
  \liminf_{|s|\nearrow+\infty} {s^{-2} F(s)} > 0 \,.
  \label{hpcoerc}
\Esist
Notice that these assumptions are fulfilled by all of the important potentials \accorpa{regpot}{obspot}.
We set, for convenience,
\Beq
  \beta := \partial\Beta , \quad
  \pi := \Pi', 
  \aand
  \Lpi := \hbox{the \Lip\ constant of $\pi$}\,.
  \label{defbetapi}  
\Eeq
Moreover, we term $D(\Beta)$ and $D(\beta)$ the effective domains of $\Beta$ and~$\beta$, respectively,
and notice that $\beta$ is a maximal monotone graph in $\erre\times\erre$.
The same symbol $\beta$ is used for the maximal monotone operators induced in $\Ldue$ and~$\LQ2$.
Finally, we introduce
\Beq
  P : H \to H, \quad
  \hbox{the orthogonal projection operator on the kernel of $B$}.
  \label{defP}
\Eeq
As for the data of our problem, we allow the forcing term appearing in \eqref{Iseconda} to 
depend on $\sigma$ and assume~that:
\Bsist
  && \fs \in \L2H;
  \label{hpfs}
  \\
  && \phiz \in \VB{\sigma_0}
  \aand
  \Beta(\phiz) \in \Luno; 
  \label{hpphiz}
  \\
  && \hbox{if $\lambda_1=0$\quad then} \quad
  \mz := \mean\phiz
  \quad \hbox{belongs to the interior of $D(\beta)$} .
  \label{hpmz}
\Esist
\Accorpa\HPdati hpfs hpmz

At this point, we make the notion of solution precise.
In the following, we use the notations
\Beq
  Q_t := \Omega \times (0,T)
  \quad \hbox{for $t\in(0,T]$}
  \aand
  Q := Q_T \,.
  \label{defQt}
\Eeq
A~solution to our system is a pair $(\phis,\mus)$ fulfilling the regularity requirements
\Bsist
  && \phis \in \H1H \cap \L\infty{\VB\sigma},
  \label{regphi}
  \\
  && \mus \in \L2{\VA{2r}},
  \label{regmu}
  \\
  && \Beta(\phis) \in \LQ1,
  \label{regBetaphi}
\Esist
\Accorpa\Regsoluz regphi regBetaphi
and satisfying the following weak formulation of the equations \accorpa{Iprima}{Icauchy}:
\Bsist
  && ( \dt\phis(t) , v )
  + ( A^r \mus(t) , A^r v )
  = 0
  \quad \hbox{for every $v\in\VA r$ and \aat},
  \qquad\quad
  \label{prima}
  \\[2mm]
  && \tau \bigl( \dt\phis(t) , \phis(t) - v \bigr)
  + \bigl( B^\sigma \phis(t) , B^\sigma( \phis(t)-v) \bigr)
  \non
  \\
  && \quad {}
  + \iO \Beta(\phis(t))
  + \bigl( \pi(\phis(t)) - \fs(t) ,  \phis(t)-v \bigr)
  \leq \bigl( \mus(t) ,  \phis(t)-v \bigr)
  + \iO \Beta(v)
  \non
  \\
  && \quad \hbox{for every $v\in\VB\sigma$ and \aat},
  \label{seconda}
  \\
  && \phis(0) = \phiz \,.
  \label{cauchy}
\Esist
\Accorpa\Pbl prima cauchy
We notice that \eqref{regBetaphi} implies that $\Beta({\phis}(t))\in\Luno$ \aat,
so that \eqref{seconda} has a precise meaning.
In the same inequality, one obviously has to read $\iO\Beta(v)=+\infty$ 
if $v\in\VB\sigma$ and $\Beta(v)\not\in\Luno$.

\Brem
\label{Moreregmu}
The regularity \eqref{regmu} of the second component of the solution is expected
even though \eqref{prima} just suggests $\mus\in\L2{\VA r}$.
Indeed, \aat\ the variational equation has the form
\Beq
  (A^r u,A^r v) = (g,v)
  \quad \hbox{for every $v\in\VA r$},
  \non
\Eeq
with $g\in H$.
From this, one easily derives that
$u\in\VA{2r}$ and $\norma{A^{2r}u}\leq\norma g$
(one can formally test by $A^{2r}u$, but a regularization procedure makes the argument rigorous).
Since $\dt\phis\in\L2H$, we thus have the regularity \eqref{regmu} as well~as
\Beq
  \dt\phis + A^{2r} \mu = 0
  \quad \aet ,
  \label{sprima}
\Eeq
i.e., the equation holds in its strong form.
\Erem

\Brem
\label{Conservation}
In the sequel, the symbol $\1$ denotes the constant function on $\Omega$ that takes the value~$1$ at every point.
With this notation, we remark that \eqref{hpsimple} implies that $A^r(\1)$ vanishes if $\lambda_1=0$,
so that \eqref{prima} {and \eqref{cauchy} imply} that
\Bsist
  && \frac d{dt} \iO \phis(t) = 0
  \quad \aat, \ \,
  \hbox{{whence}} \quad
  \non
  \\
  && \mean \phis(t) = \mz
  \quad \hbox{for every $t\in[0,T]$}
  \label{conservation}
\Esist
in this case.
On the contrary, if $\lambda_1>0$, no conservation property is expected.
\Erem

The well-posedness result (cf.~\cite[Thm.~2.6]{CGS18}) reads as follows:

\Bthm
\label{Wellposedness}
Let the assumptions \eqref{hpAB}, \accorpa{hprst}{hpVB} and 
\accorpa{hpF}{hpcoerc}
on the structure of the system
and \HPdati\ on the data be fulfilled.
Then, there exists a pair $(\phis,\mus)$ satisfying \Regsoluz\
and solving problem \Pbl.
Moreover, the component $\phis$ of the solution is unique.
\Ethm

\Brem
\label{Uniquemu}
No uniqueness for the component $\mus$ of the solution can be expected, in general.
However, in particular situations, $\mus$~is unique, too.
This is the case if $\lambda_1>0$.
Indeed, this assumption implies that $A^{2r}$ is invertible so that 
\eqref{sprima} can be uniquely solved for~$\mus$.
On the contrary, the case $\lambda_1=0$ is much more delicate.
A~sufficient condition that ensures uniqueness for $\mus$ is the following (see \cite[Rem.~4.1]{CGS18}):
$\Beta$~is an everywhere defined $C^1$ function and $\phis$ is bounded.
We notice that the same argument used in the quoted remark also applies if $D(\beta)$ is an open interval and 
$\beta$ is a continuous single-valued function on~it 
(like in the case \eqref{logpot} of the logarithmic potential)
provided that all {of} the values of $\phis$ belong to a compact subset of~$D(\beta)$.
\Erem

Let us come to the results of this paper.
The first deals with the \bhv\ of the solutions to problem \Pbl\ as $\sigma$ \newjuerg{tends} to zero.

\Bthm
\label{Asymptotics}
Besides the assumptions of Theorem~\ref{Wellposedness},
assume that
\Beq
  \fs \to f 
  \quad \hbox{strongly in $\L2H$ as $\sigma\searrow0$}.
  \label{hpf}
\Eeq
Then, for every $\sigma>0$ there is a 
solution $(\phis,\mus)$ to problem \Pbl\ such that
\Bsist
  && \phis \to \phi
  \quad \hbox{weakly in $\H1H$},
  \label{convphi}
  \\
  && \mus \to \mu
  \quad \hbox{weakly in $\L2{\VA{2r}}$},
  \label{convmu}
  \\
  && B^\sigma\phis \to \zeta
  \quad \hbox{weakly star in $\L\infty H$},
  \label{convBs}
\Esist
\Accorpa\Convsoluz convphi convBs
as $\sigma\searrow0$, possibly along a subsequence, for some triplet $(\phi,\mu,\zeta)$ satisfying
\Beq
  \phi \in \H1H , \quad
  \mu \in \L2{\VA{2r}} 
  \aand
  \zeta\in \L\infty H \,.
  \label{regsoluzz}
\Eeq
Moreover, under the additional assumption,
\begin{align}
&\hbox{for all } v\in H  \hbox{ such that } \Beta(v) \in L^1(\Omega)
\hbox{ there exists a sequence } \{v_n\} \subset V_B^{\sigma_0} 
\non
\\
&\hbox{such that } \
v_n \to v \ \hbox{ in } H \ \hbox{ and } \ \liminf_{n\to \infty } \iO \Beta(v_n) = \iO \Beta(v)\,, 
\label{ass-pier}
\end{align}
\newjuerg{the following holds true:} whenever $(\phis,\mus)$ is a solution to problem \Pbl\ for $\sigma>0$ 
and \Convsoluz\ hold for some {triplet} $(\phi,\mu,\zeta)$,
then $\zeta=\phi-P\phi$ and the pair $(\phi,\mu)$ is a solution to the system
\Bsist
  && ( \dt\phi(t) , v )
  + ( A^r \mu(t) , A^r v )
  = 0
  \quad \hbox{for every $v\in\VA r$ and \aat},
  \qquad\quad
  \label{primaz}
  \\[2mm]
  && \tau \bigl( \dt\phi(t) , \phi(t) - v \bigr)
  + \bigl( \phi(t) - P \phi(t) , \phi(t)-v \bigr)
  \non
  \\
  && \quad {}
  + \iO \Beta(\phi(t))
  + \bigl( \pi(\phi(t)) - f(t) ,  \phi(t)-v \bigr)
  \leq \bigl( \mu(t) ,  \phi(t)-v \bigr)
  + \iO \Beta(v)
  \non
  \\
  && \quad \hbox{for every $v\in H$\, and \aat},
  \label{secondaz}
  \\
  && \phi(0) = \phiz \,.
  \label{cauchyz}
\Esist
\Accorpa\Pblz primaz cauchyz
\Ethm

\Brem
\label{Remasym}
The above statement looks a little involved. 
Besides the assumption~\eqref{ass-pier} 
we are going to discuss in a while, we point out that 
that no uniqueness for the solution $(\phis,\mus)$ is required.
On the contrary, if additional assumptions were made that guarantee uniqueness for $(\phis,\mus)$
(see Remark~\ref{Uniquemu}) and \eqref{ass-pier} \newjuerg{were} assumed, 
then the statement would look much simpler, namely: 
as $\sigma$ tends to zero, the solution $(\phis,\mus)$ converges 
(in the sense of \accorpa{convphi}{convmu}, possibly along a subsequence)
to~a solution $(\phi,\mu)$ to problem \Pblz.
If, in addition, uniqueness holds for the solution $(\phi,\mu)$ to the limiting problem,
then the whole family $\{(\phis,\mus)\}$ converges to $(\phi,\mu)$ as $\sigma$ tends to zero.
\Erem

\Brem
\label{No-ass-pier}
As observed in the forthcoming Remark~\ref{No-ass-pier-bis},
if \eqref{ass-pier} is not assumed, a weaker conclusion can anyway be obtained:
the variational inequality \eqref{secondaz} is fulfilled by all the test functions $v\in\VB{\sigma_0}$.
Indeed, it is stressed in the remark that 
assumption \eqref{ass-pier} is used in the proof of Theorem~\ref{Asymptotics} just to extend to any $v\in H$ 
the validity of \eqref{secondaz} proved for test functions $v\in\VB{\sigma_0}$.
\Erem

\Brem
\label{From-ass-pier}
So, if \eqref{ass-pier} is assumed, \newjuerg{then}
every limiting pair $(\phi,\mu)$ satisfies \eqref{secondaz} with arbitrary test functions $v\in H$.
This has the following important consequence:
there exists \newjuerg{some} $\xi$ satisfying
\Bsist
  && \xi \in \L2H
  \aand
  \xi \in \beta(\phi)
  \quad \aeQ \,,
  \label{regxi}
  \\
  && \tau \dt\phi + \phi - P\phi + \xi + \pi(\phi) = \mu + f
  \quad \aeQ \,.
  \label{eqsecondaz}  
\Esist
Indeed, if we set
\Beq
  \xi := \mu - \tau \dt\phi - \phi + P\phi - \pi(\phi) + f
  \label{defxi}
\Eeq
then $\xi$ belongs to $\L2H$, equation \eqref{eqsecondaz} is satisfied, and \eqref{secondaz} becomes
\Beq
  \iO \Beta(\phi(t)
  \leq \bigl( \xi, \phi(t) - v \bigr)
  + \iO \Beta(v)
  \quad \hbox{for every $v\in H$\, and \aat}.
  \label{xisubdiff}
\Eeq
But this exactly means that $\xi(t)\in\partial\Beta(\phi(t))=\beta(\phi(t))$ \aat, i.e., the second condition 
in~\eqref{regxi}.
If instead \eqref{ass-pier} is not assumed, then \eqref{secondaz} is satisfied only for test functions $v\in\VB{\sigma_0}$,
as said in Remark~\ref{No-ass-pier}.
Nevertheless, the definition \eqref{defxi} still yields $\xi\in\L2H$ and implies that \eqref{eqsecondaz} is satisfied.
However, in this case, \eqref{xisubdiff} is only true for $v\in\VB{\sigma_0}$,
and this means that \aat\ the function $\xi(t)$ belongs to the subdifferential 
of the function $\VB{\sigma_0}\ni v\mapsto\iO\Beta(v)$.
Notice that this subdifferential is a subset of the dual space~$(\VB{\sigma_0})^*$
and might contain elements that do not belong to~$H$
(in~the sense of the Hilbert triplet $(\VB{\sigma_0},H,(\VB{\sigma_0})^*)$).
Moreover, if a function $u\in H$ belongs to such a subdifferential, \newjuerg{then}
it is not clear whether it also belongs to the subdifferential in~$H$
(i.e., that of the function $H\ni v\mapsto\iO\Beta(v)$),
so that we cannot conclude that $\xi\in\beta(\phi)$ \aeQ. \pier{About this matter, let us quote the paper \cite{BCGG} for related issues.}

\Erem

\Brem
\label{Suff-ass-pier}
A sufficient condition for \eqref{ass-pier} \newjuerg{to hold true} 
is~the following
(satisfied in all of the concrete cases, at least if $\sigma_0$ \newjuerg{is} small enough):
\Beq
  \Hdue \subset \VB{\sigma_0} \,.
  \label{suff-ass-pier}
\Eeq
In order to construct the sequence $\{\vn\}$ for a given $v\in H$, 
we solve the Neumann boundary value problem
\Beq
  \iO \vn z + \frac 1n \iO \nabla\vn \cdot \nabla z
  = \iO vz
  \quad \hbox{for every $z\in\Huno$}.
  \label{singpert}
\Eeq
Since $v\in H$, we have that $\vn\in\Hdue$ and thus $\vn\in\VB{\sigma_0}$, by~\eqref{suff-ass-pier}.
\newjuerg{Now, if we take $z=\vn$ in \eqref{singpert} \pier{and use the Cauchy--Schwarz inequality in the \rhs}, then we easily find that}
\Beq
\newjuerg{\norma\vn\leq\norma v \quad\mbox{and}\quad 
\norma{\mbox{$\frac 1n$}\nabla\vn}^{\pier 2}\leq 
\pier{\textstyle \frac 1{n}} \norma{v}^{\pier 2}
\quad \mbox{for all $n\in\enne$.}}
\label{paolauno}
\Eeq
\newjuerg{Hence, there are a subsequence $\{v_{n_k}\}$ and some $w\in H$ such that $v_{n_k}\to w$ weakly in $H$.
Moreover, since $\frac 1n \nabla\vn\to( 0,0,0)$ strongly in $H\times H\times H$ by \eqref{paolauno}, we easily infer from 
\eqref{singpert} that $w=v$. A fortiori, by the uniqueness of the limit point, the entire sequence $\{v_n\}$ 
converges weakly in $H$ to $v$. But then, by the weak sequential lower semicontinuity of norms,} 
$$\newjuerg{\|v\|\le \liminf_{n\to\infty}\|\vn\|\le\limsup_{n\to\infty}\|\vn\|\le\|v\|,}$$
\newjuerg{where the latter inequality follows from \eqref{paolauno}. We thus have 
\,$\|v\|=\lim_{n\to\infty}\|\vn\|$, and the uniform convexity of $H$ yields that $\vn\to v$ strongly in $H$.}

Now, denoting by $\Betaeps$ and $\betaeps$ the Moreau--Yosida $\eps$-approximations 
of $\Beta$ and~$\beta$, respectively 
(see, e.g., \cite[p.~28]{Brezis}),
we~account for the definition of the subdifferential $\betaeps=\partial\Betaeps$
and the identity obtained by testing \eqref{singpert} by $\betaeps(\vn)\in\Huno$. 
We have~that
\Beq
  \iO \Betaeps(\vn) - \iO \Betaeps(v)
  \leq \iO \betaeps(\vn) (\vn-v)
  = - \frac 1n \iO \betaeps'(\vn) |\nabla\vn|^2
  \leq 0
  \non
\Eeq
and we deduce that
\Beq
  \iO \Betaeps(\vn)
  \leq \iO \Betaeps(v)
  \leq \iO \Beta(v) 
  \quad \hbox{whence also} \quad
  \iO \Beta(\vn)
  \leq \iO \Beta(v) 
  \non
\Eeq
by letting $\eps$ tend to zero.
This implies the inequality ``$\leq$'' in \eqref{ass-pier}.
Since the opposite inequality {clearly follows from the
lower semicontinuity of the function $z\mapsto \iO \Beta (z)$ in $H$, we finally conclude for the validity of  \eqref{ass-pier}.}
\Erem

Notice that Theorem~\ref{Wellposedness} ensures the existence of at least one solution to the limiting problem \Pblz\
with the regularity specified in~\eqref{regsoluzz}.
Our next result deals with partial uniqueness and continuous dependence of the solution.
This {will be}  proved in the last Section~\ref{LIMITING},
{which is devoted to the study of the limiting problem}.

\Bthm
\label{Partuniqz}
Let the general assumptions on the structure be fulfilled,
and assume that $\phiz$ satisfies \eqref{hpphiz}.
Morever, let $f_i\in\L2H$, $i=1,2$, be two choices  of the forcing term $f$ appearing in \eqref{secondaz},
and let $(\phi_i,\mu_i)\in\H1H\times\L2{\VA{2r}}$ be two corresponding solutions to problem {\Pblz}\ with $f=f_i$. 
Then we~have
\Beq
  \norma{\phi_1-\phi_2}_{\L\infty H}
  \leq C_{cd} \norma{f_1-f_2}_{\L2H},
  \label{contdep}
\Eeq
with a constant $C_{cd}$ that depends only on~$\tau$, the \Lip\ constant~$\Lpi$, and~$T$.
In particular, if $f\in\L2H$, the first component $\phi$ of the solution $(\phi,\mu)$ to problem \Pblz\ 
is uniquely determined.
\Ethm

In our final result, we require some regularity of the data and further assumptions on the structure
that are satisfied in all of the concrete cases of interest,
and we prove a regularity result.
As a byproduct, we obtain a sufficient condition for the uniqueness of the second component $\mu$ of the solution.
Sufficient conditions for uniqueness in a different direction are given in the forthcoming Remark~\ref{Uniquemuz}.

\Bthm
\label{Eqsecondaz} 
Let the general assumptions on the structure be fulfilled.
In addition, assume {that}
\begin{align}
  & \VB n \subset \Huno
  \quad \hbox{for some positive integer $n$,}
  \label{newhpsuppl}
  \\[2mm]
  & \VA{2r} \subset \Hx\eta \,, \quad
  f \in \L2{\Hx\eta}
  \aand
  \phiz \in \Hx\eta
  \quad \hbox{for some $\eta\in(0,1]$,}
  \label{hpsuppl}
\end{align}
and let $(\phi,\mu)$ {with}
\Beq
  \phi \in \H1H
  \aand
  \mu \in \L2{\VA{2r}}
  \label{regsoluz}
\Eeq
be a solution to problem \Pblz.
Then $\phi$ enjoys the further regularity
\Beq
  \phi \in \L2{\Hx\eta},
  \label{moreregphi}
\Eeq
{and there exists some $\xi$ satisfying \accorpa{regxi}{eqsecondaz}.}
In particular, even the second component $\mu$ of the solution is unique if $\beta$ is single-valued.
\Ethm

Throughout the paper, we widely use the Cauchy--Schwarz and Young inequalities,
the latter in the form
\Beq
  ab\leq \delta a^2 + \frac 1 {4\delta}\,b^2
  \quad \hbox{for every $a,b\in\erre$ and $\delta>0$}.
  \label{young}
\Eeq
Moreover, in performing a priori estimates, we use the same small letter $c$
for (possibly) different constants that depend only 
on the structure of our system but~$\sigma$, and on the assumptions on the data.
In particular, the values of $c$ do not depend on the regularization parameter $\lambda$ we introduce in the next section.
On the contrary, some precise constants are denoted by different symbols
(see, e.g., \eqref{poincare}, where a capital letter with an index is used).


\section{Asymptotic analysis}
\label{ASYMPTOTICS}
\setcounter{equation}{0}

This section is devoted to the proof of Theorem~\ref{Asymptotics}.
The construction of the solutions $(\phis,\mus)$ mentioned in the statement
relies on a priori estimates on the solutions to a regularized problem,
as done in~\cite{CGS18} to solve problem \Pbl\ with a fixed~$\sigma$.
Hence, we briefly recall that regularization procedure.
For $\lambda>0$ (small enough if needed), let $\betal$ be the Yosida approximation of $\beta$ at the level~$\lambda$
(see, e.g., \cite[p.~28]{Brezis}).
The corresponding Moreau regularization $\Betal$ of $\Beta$ is thus given~by
\Beq
  \Betal(s) = \int_0^s \betal(s') \, ds'
  \quad \hbox{for $s\in\erre$},
  \non
\Eeq
since $\,\betal(0)=0\,$ due to \eqref{hpBeta}.
Then, the regularized problem consists in looking for a pair $(\phisl,\musl)$
satisfying the regularity requirements 
\Beq
  \phisl \in \H1H \cap \L\infty{\VB\sigma} \cap \L2{\VB{2\sigma}}
  \aand
  \musl \in \L2{\VA{2r}},
  \label{regsoluzl}
\Eeq
and solving the following system:
\Bsist
  && ( \dt\phisl(t) , v )
  + ( A^r \musl(t) , A^r v )
  = 0
  \quad \hbox{for every $v\in\VA r$ and \aat},
  \qquad
  \label{primal}
  \\[1mm]
  \separa
  && {\tau\bigl(\dt\phisl(t),v\bigr)}
  + \bigl( B^\sigma\phisl(t) , B^\sigma v \bigr)
  + \bigl( \betal(\phisl(t)) + \pi(\phisl(t)) - \fs(t) ,  v \bigr)
  = \bigl( \musl(t) ,  v \bigr)
  \non
  \\[1mm]
  && \quad \hbox{for every $v\in\VB\sigma$ and \aat}{,}
  \label{secondal}
  \\[1mm]
  && \phisl(0) = \phiz \,.
  \label{cauchyl}
\Esist
\Accorpa\Pbll primal cauchyl
We notice that the variational inequality \eqref{seconda}
is replaced by the equality \eqref{secondal} in the approximating problem
(since $\betal$ is an everywhere defined \Lip\ continuous function).
{The existence part of Theorem~\ref{Wellposedness}
is proved by solving the above regularized problem
(cf. \cite[Thm.~5.1]{CGS18})}
and showing that its solution $(\phisl,\musl)$ converges as $\lambda\searrow0$
(in~a suitable topology, possibly just along a subsequence)
to~a pair $(\phis,\mus)$ which turns out to solve problem \Pbl.
This solution, where now $\sigma$ is a varying parameter that we \newjuerg{intend to approach zero},
will be the good candidate for Theorem~\ref{Asymptotics}.

Before starting {to estimate}, it is worth observing that Remark~\ref{Moreregmu}
applies to both equations \eqref{primal} and \eqref{secondal}.
This is obvious for the former.
As {far as} the latter is concerned, one has to replace $A$ and $r$ by $B$ and~$\sigma$, respectively,
and notice that $\betal$ is \Lip\ continuous, so that
$\musl+\fs-\betal(\phisl)-\pi(\phisl)\in\L2H$.
This justifies the last regularity condition for $\phisl$ in \eqref{regsoluzl}
(in~contrast with \eqref{regphi}) 
and implies the strong form of both equations,~i.e.,
\Bsist
  && \dt\phisl
  + A^{2r} {\musl} = 0
  \quad \aet,
  \label{sprimal}
  \\
  && \tau \dt\phisl
  + B^{2\sigma} \phisl
  + \betal(\phisl)
  + \pi(\phisl)
  = \musl
  + \fs
  \quad \aet \,.
  \label{ssecondal}
\Esist
We also recall the convention on the symbol $c$ for possibly different constants
made at the end of Section~\ref{STATEMENT}.
Moreover, since \eqref{hpf} implies that $\fs$ is bounded in $\L2H$,
we allow $c$ to also depend on a bound for the corresponding norm.

\step
First a priori estimate

We test \eqref{primal} written at the time~$s$ by~$\mus(s)$.
At the same time, we {insert $+ \phisl (s) $ to both sides of 
\eqref{ssecondal} written at the time~$s$ and multiply it by~$\dt\phisl(s)$, then integrating over~$\Omega$.
We sum up both equalities,} noting that a cancellation occurs, and integrate over $(0,t)$ with respect to~$s$.
We obtain
\Bsist
  && \iot \norma{A^r \musl(s)}^2 \, ds
  + \tau \intQt |\dt\phisl|^2
  + \frac 12 \, { \bigl(\norma{\phisl(t)}^2 + \norma{B^\sigma \phisl(t)}^2\bigr)}
  + \iO \Betal(\phisl(t))
  \non
  \\
  && {} = \frac 12 \, { \bigl(\norma{\phiz}^2 + \norma{B^\sigma \phiz}^2\bigr)}
  + \iO \Betal(\phiz) 
  + \intQt (\fs {{}+ \phisl}- \pi(\phisl)) \dt\phisl \,.
  \non
\Esist
Even the last integral on the \lhs\ is nonnegative.
We estimate the terms on the \rhs\ 
by accounting for {the} assumptions \eqref{hpphiz} and \eqref{hpf} on $\phiz$ and~$\fs$, respectively,
and owing to the \Lip\ continuity of~$\pi$. 
{Recalling also \eqref{defprodBs},} we have that
\begin{align}
  &{ \norma{\phiz}^2 + \norma{B^\sigma \phiz}^2 = \norma\phiz_{B,\sigma}^2
  ={}} \somma j1\infty (1+ (\lambda'_j)^{2\sigma}) |(\phiz,e'_j)|^2
  \non
  \\[-2mm]
  & \leq \somma j1\infty ({2} + (\lambda'_j)^{2\sigma_0}) |(\phiz,e'_j)|^2\leq {2} \norma\phiz_{B,\sigma_0}^2\,,
  \non
  \\[2mm]
  \separa
  & \iO \Betal(\phiz) \leq \iO \Beta(\phiz)\,, 
  \non
  \\[2mm]
  \separa
  & \intQt (\fs {{}+ \phisl} - \pi(\phisl)) \dt\phisl
  \leq \frac \tau 2 \intQt |\dt\phisl|^2
  + c \iot \bigl({ \|\fs(s)\|^2 + \|\phisl(s)\|^2} + 1 \bigr) \, ds
  \non
  \\
  & \leq \frac \tau 2 \intQt |\dt\phisl|^2
  + c \iot {\|\phisl(s)\|^2} \, ds + c \,. \non
\end{align}
Therefore, by rearranging and applying the Gronwall lemma, we conclude~that
\Beq
  \norma{A^r\musl}_{\L2H}
  + \norma\phisl_{\H1H}
  + {\norma{\phisl}_{\L\infty {V_B^\sigma}}}
  + \norma{\Betal(\phisl)}_{\L\infty\Luno}
  \leq c \,.
  \label{primastima}
\Eeq
From this and \eqref{sprimal} we deduce that
\Beq
  \norma{A^{2r}\musl}_{\L2H}
  \leq c \,.
  \label{daprimastima}
\Eeq

\step
Second a priori estimate

Our aim is {to improve} the estimate concerning~$\musl$.
Indeed, for the {following}     we need that
\Beq
  \norma\musl_{\L2{\VA r}} \leq c \,.
  \label{secondastima}
\Eeq
We notice at once that this and \eqref{daprimastima} would imply that
\Beq
  \norma\musl_{\L2{\VA{2r}}} \leq c \,.
  \label{dasecondastima}
\Eeq
The desired estimate trivially follows from \eqref{primastima} if $\lambda_1>0$.
So, we now deal with the {other} case $\lambda_1=0$
and apply a well-known trick
based on the assumption \eqref{hpmz} and the consequent inequality
\Beq
  \betal(s) (s-\mz)
  \geq \delta_0 |\betal(s)| - C_0, 
  \label{trickMZ}
\Eeq
which holds for some $C_0>0$ and every $s\in\erre$ and $\lambda\in(0,1)$,
where $\delta_0$ is such that the interval $[\mz-\delta_0,\mz+\delta_0]$ 
is included in the interior of~$D(\beta)$
(cf.\ \cite[Appendix, Prop. A.1]{MiZe}; see also \cite[p.~908]{GiMiSchi} for a detailed proof).
Inequality \eqref{trickMZ} implies that
\Beq
  \bigl( \betal(\phisl(t)) , \phisl(t)-\mz{\1} \bigr)
  \geq \delta_0 \, \norma{\betal(\phisl(t))}_{\Luno} - c 
  \quad \aat ,
  \label{fromMZ}
\Eeq
and this can be used when testing equation \eqref{secondal} by ${\phisl}-\mz{\1}$.
To this concern, we recall that $\1\in\VB\sigma$ by \eqref{hpVB}
and notice that the conservation property \eqref{conservation}
also holds for~$\phisl$, i.e., $\mean\phisl(t)=\mz$ for every $t\in[0,T]$.
So, \aat, we test \eqref{secondal} by $\phisl(t)-\mz{\1}$ and rearrange a little.
However, we omit writing the time $t$ for a while.
We also write $k$ instead of $k\1$ if $k$ is a real number.
We have \aet\ that
\Bsist
  &&\norma{B^\sigma\phisl}^2
  + \bigl( \betal(\phisl) , \phisl-\mz \bigl)
  \non
  \\
  && = ( \musl , \phisl-\mz )
  + \bigl( \fs - \tau\dt\phisl - \pi(\phisl) , \phisl-\mz \bigr) 
  + (B^\sigma\phisl , B^\sigma\mz).
  \label{persecondastima}
\Esist
The \lhs\ of this equality can be estimated from below by {virtue of} \eqref{fromMZ}.
The first term on the \rhs\ can be dealt with by accounting for the Poincar\'e type inequality \eqref{poincare} 
{as follows:}
\Bsist
  && ( \musl , \phisl-\mz )
  = ( \musl - \mean\musl \,, \phisl-\mz )
  \leq \norma{\musl-\mean\musl} \, \norma{\phisl-\mz}
  \non
  \\
  && \leq c \, \norma{A^r(\musl-\mean\musl)} \, \norma{\phisl-\mz}
  = c \, \norma{A^r\musl} \, \norma{\phisl-\mz}\,, 
  \non
\Esist
the last equality being due to $A^r\1=0$.
Therefore, by recalling \eqref{primastima}, we have that 
the whole \rhs\ of \eqref{persecondastima} is bounded in~$L^2(0,T)$
and conclude~that
\Beq
  \norma{\betal(\phisl)}_{\L2\Luno} \leq c, 
  \quad \hbox{whence immediately} \quad
  \norma{\mean\betal(\phisl)}_{L^2(0,T)} \leq c \,.
  \non
\Eeq
At this point, we can test the second equation \eqref{secondal} by $\1$
and deduce a bound for $\mean\musl$ in~$L^2(0,T)$.
This and \eqref{primastima} imply~\eqref{secondastima}.
As already noticed, \eqref{dasecondastima} is proved as well.

\step
First conclusion

As already remarked, in the proof of \cite[Thm.~5.1]{CGS18} with a fixed $\sigma$ it is shown that
$(\phisl,\musl)$ converges as $\lambda$ tends to zero 
(in a proper topology, possibly along a subsequence)
to~some pair $(\phis,\mus)$, 
and it is proved that such a pair is a solution to problem \Pbl.
We prove that the family $\{(\phis,\mus)\}_{\sigma>0}$ constructed \newjuerg{in} this way
satisfies all the requirement of the statement.
The starting point is the conservation of the bounds just proved in the limit as $\lambda\searrow0$.
We have~that
\Beq
  \norma\phis_{\H1H}
  + \norma\mus_{\L2{\VA{2r}}}
  + \norma{B^\sigma\phis}_{\L\infty H}
  \leq c \,,
  \non
\Eeq
and we conclude that \Convsoluz\ hold true for some triplet $(\phi,\mu,\zeta)$ satisfying \eqref{regsoluzz}.
This ends the proof of the first part of the statement.

Let us come to the second part.
So, we assume that $\{(\phis,\mus)\}_{\sigma>0}$ is a family of solutions to problem \Pbl\
and that \Convsoluz\ hold true for some triplet $(\phi,\mu,\zeta)$ satisfying \eqref{regsoluzz} as $\sigma\searrow0$,
possibly for a subsequence
(however, we always write $\sigma$ instead of the elements of some subsequence $\{\sigma_k\}$, for brevity).
We have to prove that $\zeta=\phi-P\phi$ and that $(\phi,\mu)$ solves problem~\Pblz\pier{, by also assuming \eqref{ass-pier}.}

\step
First characterization

We {are going to} show that $\zeta=\phi-P\phi$ by proving that
\Beq
  B^\sigma \phis \to \phi-P\phi
  \quad \hbox{weakly in $\LQ2$}.
  \label{basic}
\Eeq
{To this end, we} use the eigenvalues $\lambda'_j$ and the eigenfunctions $e'_j$ of $B$
and notice that $e'_j$ is orthogonal to $\ker B$ if $\lambda'_j>0$
while $\lambda'_j=0$ if $e'_j\in\ker B$.
We set, for convenience,
\Beq
  \Asj := \ioT \bigl( B^\sigma \phis(t) , \psi(t) \, e'_j \bigr) \, dt
  = (\lambda'_j)^\sigma\ioT \bigl( \phis(t) , \psi(t) \, e'_j \bigr) \, dt
  \non
\Eeq
for $\psi\in L^2(0,T)$ and $j=1,2,\dots$,
and we notice that \eqref{basic} follows if we prove that
\Beq
 \lims \Asj
  = \Azj
  := \ioT \bigl( \phi(t) - P\phi(t) , \psi(t) \, e'_j \bigr) \, dt
  \label{charact}
\Eeq
for every $\psi$ and~$j$ as before,
since the linear combinations of the products $\psi\,e'_j$ of such real functions and eigenfunctions of $B$
form a dense subspace of~$\LQ2$.
So, we fix $\psi$ and~$j$.
As for~$j$, we distinguish two cases.
Assume first that $\lambda'_j>0$.
Then, $(\lambda'_j)^\sigma$ tends to $1$ as $\sigma$ tends to zero.
Moreover, \eqref{convphi} holds.
We thus deduce that
\Beq
 \lims \Asj
  = \ioT \bigl( \phi(t) , \psi(t) \, e'_j \bigr) \, dt
  = \Azj ,
  \non
\Eeq
the last equality being due to the orthogonality between $P\phi(t)$ and~$e'_j$.
Assume now that $\lambda_j=0$.
Then, we trivially have that $\Asj=0$ for every $\sigma>0$.
On the other hand, we also have that $\Azj=0$ since 
$e'_j\in\ker B$ and $\phi(t)-P\phi(t)$ is orthogonal to $\ker B$ \aat.
Therefore, \eqref{charact} is proved in any case.

\Brem
\label{ConvBsv}
The same argument shows that, for every $v\in\L2{\VB{\sigma_0}}$,
the weak convergence $B^\sigma v\to v-Pv$ in $\L2H$ holds true as $\sigma$ tends to zero.
In fact, the convergence is strong:
\Beq
  B^\sigma v \to v-Pv
  \quad \hbox{strongly in $\L2H$ \ for every $v\in\L2{\VB{\sigma_0}}$} .
  \label{convBsv}
\Eeq
Indeed, \aat, $B^\sigma v(t) \to v(t)-Pv(t)$ strongly in $H$ by \cite[Lem.~7.5]{CG1}.
Moreover, the Lebesgue dominated convergence theorem can be applied since
\Beq
  \norma{B^\sigma v(t)}^2
  = \somma j1\infty (\lambda'_j)^{2\sigma} |(v(t),e'_j)|^2
  \leq \somma j1\infty (1 + (\lambda'_j)^{2\sigma_0}) |(v(t),e'_j)|^2
  = \norma{v(t)}_{B,{\sigma_0}}^2
  \non
\Eeq
\aat\ and every $\sigma\in(0,\sigma_0]$, and $\norma{v(\cpto)}_{B,{\sigma_0}}^2$ belongs to $L^1(0,T)$.
\Erem

\medskip

To conclude the proof, we have to show that $(\phi,\mu)$ solves problem \Pblz\
{under the further assumption~\eqref{ass-pier}}.
The first equation obviously follows from~\eqref{prima} due to \accorpa{convphi}{convmu},
and the initial condition \eqref{cauchyz} is satisfie{d} as well 
since \eqref{convphi} implies weak convergence in $\C0H$.
So, it remains to verify the variational inequality \eqref{secondaz}.
To this concern, it is convenient to give different formulations of both \eqref{seconda} and~\eqref{secondaz}.
This {procedure} is based on the lemma stated below,
which follows from the classical theory of variational inequalities of elliptic type
in the framework of Convex Analysis.
However, for the reader's convenience, we also give {a} simple proof.

\Blem
\label{Known}
Let $V$ be a Hilbert space, $\Vp$ its dual space, $\<\cpto,\cpto>$ the duality pairing between $\Vp$ and~$V$, 
and $a:V\times V\to\erre$ a continuous bilinear form.
Moreover, assume that
\Bsist
  && \Gam_1 : V \to (-\infty,+\infty] 
  \quad \hbox{is convex, proper and lower semicontinuous},
  \qquad
  \label{hpG1}
  \\
  && \Gam_2 : V \to \erre
  \quad \hbox{is convex and G\^ateaux differentiable},
  \non
  \\
  && \quad
  \hbox{and \ $\gam_2:V\to\Vp$ \ is its G\^ateaux derivative.}
  \label{hpG2}
\Esist
Then, for every $u\in V$ and $g\in\Vp$, the variational inequalities
\Bsist
  && a(u,u-v) + \Gam_1(u) + \< \gam_2(u),u-v >
  \leq \< g,u-v> + \Gam_1(v)
  \quad \hbox{for every $v\in V$},
  \qquad
  \label{oldvi}
  \\
  && a(u,u-v) + \Gam_1(u) + \Gam_2(u)
  \leq \< g,u-v> + \Gam_1(v) + \Gam_2(v)
  \quad \hbox{for every $v\in V$},
  \label{newvi}
\Esist
are equivalent to each other.
\Elem

\Bdim
Assume \eqref{oldvi} and let $v\in V$.
Since $\Gam_2$ is convex and $\gam_2$ is its derivative, we have~that
\Beq
  \Gam_2(u) \leq \< \gam_2(u),u-v > + \Gam_2(v),
  \non
\Eeq
whence the chain
\Bsist
  && a(u,u-v) + \Gam_1(u) + \Gam_2(u)
  \non
  \\
  && \leq a(u,u-v) + \Gam_1(u) + \< \gam_2(u) , u-v > + \Gam_2(v)
  \non
  \\  
  && \leq \< g , u-v > + \Gam_1(v) + \Gam_2(v) 
  \non
\Esist
{follows,} that is, \eqref{newvi}.
Assume now \eqref{newvi} and let $v\in V$.
By writing \eqref{newvi} with $w$ in place of $v$ 
and then choosing $w=u+\theta(v-u)$ with $\theta\in(0,1)$ (whence $u-w=\theta(u-v)$),
we obtain~that
\Beq
  \theta \, a(u,u-v) + \Gam_1(u) + \Gam_2(u)
  \leq \theta \, \< g,u-v> + \Gam_1(u+\theta(v-u)) + \Gam_2(u+\theta(v-u)).
  \non
\Eeq
By rearranging and dividing by~$\theta$, we deduce that
\Beq
  a(u,u-v)
  + \frac {\Gam_1(u) - \Gam_1(u+\theta(v-u))} \theta
  + \frac {\Gam_2(u) - \Gam_2(u+\theta(v-u))} \theta
  \leq \< g,u-v> \,.
  \non
\Eeq
On the other hand, the convexity of $\Gam_1$ implies that
\Beq
  \Gam_1(u)
  \leq \frac {\Gam_1(u) - \Gam_1(u+\theta(v-u))} \theta
  + \Gam_1(v) \,.
  \non
\Eeq
By combining these inequalities, we deduce that
\Beq
  a(u,u-v)
  + \Gam_1(u) 
  + \frac {\Gam_2(u) - \Gam_2(u+\theta(v-u))} \theta
  \leq \< g,u-v >
  + \Gam_1(v) ,
  \non
\Eeq
and letting $\theta$ tend to zero, we obtain~\eqref{oldvi}.
\Edim

As already announced, we use this lemma to replace both \eqref{seconda} and \eqref{secondaz}
by different variational inequalities.

\step
{First alternative formulation}

We first observe that \eqref{seconda} for every $v\in\VB\sigma$ as required
implies the same inequality for every $v\in\VB{\sigma_0}$ 
since $\VB{\sigma_0}\subset\VB\sigma$.
Now, by recalling that $\Lpi$ is the \Lip\ constant of~$\pi$,
we replace the latter variational inequality by an equivalent one 
by applying lemma with the choices 
\Bsist
  & V = \VB{\sigma_0} \,, \quad 
  a(u,v) = \iO ( B^\sigma u,B^\sigma v) 
  - \Lpi (u,v) 
  \quad \hbox{for $u,v\in V$},
  \non
  \\
  & \Gam_1(v) = \iO \Beta(v)
  \aand
  \Gam_2(v) = \iO \bigl( \Pi(v) + \frac \Lpi 2 \, v^2 \bigr)
  \quad \hbox{for $v\in V$},
  \non
  \\
  & \quad
  \hbox{and, \aat,} \quad
  u = \phis(t)
  \aand
  g = \mus(t) + \fs(t) - \tau \dt\phis(t) \,.
  \non
\Esist
Notice that $\Gam_2$ actually is convex (since $\pi'+\Lpi\geq0$ {a.e. in $\erre$}) and G\^ateaux differentiable 
and that its derivative $\gam_2$ is given by
$\<\gam_2(u),v>=(\pi(u)+\Lpi u,v)$.
Hence, we deduce that the variational inequality \eqref{seconda} 
{required just for every $v\in\VB{\sigma_0}$} is equivalent~to
\Bsist
  && \tau \bigl( \dt\phis(t) , \phis(t) - v \bigr)
  + \bigl( B^\sigma \phis(t) , B^\sigma( \phis(t)-v) \bigr)
  - \Lpi \bigl( \phis(t) , \phis(t)-v \bigr)
  \non
  \\
  && \quad {}
  + \iO \Alpha(\phis(t))
  \leq \bigl( \mus(t) + \fs(t) ,  \phis(t)-v \bigr)
  + \iO \Alpha(v)
  \non
  \\
  && \quad \hbox{for every $v\in\VB{\sigma_0}$ and \aat} ,
  \label{newseconda}
\Esist
where, for brevity, we have set 
\Bsist
  \Alpha(s) := \Beta(s) + \Pi(s) + \frac \Lpi 2 \, s^2
  \quad \hbox{for $s\in\erre$} \,.
  \label{defAlpha}
\Esist
{We fix what we have established:
\Beq
  \hbox{\sl the variational inequality \eqref{seconda} implies \eqref{newseconda}\/.}
  \label{replseconda}
\Eeq
\vskip -6pt
}%

\step
{Second alternative formulation}

{Similarly, we would like to show that \eqref{secondaz} is equivalent~to}
\Bsist
  && \tau \bigl( \dt\phi(t) , \phi(t) - v \bigr)
  + \bigl( \phi(t) - P \phi(t) , \phi(t)-v \bigr)
  - \Lpi \bigl( \phi(t) , \phi(t)-v \bigr)
  \non
  \\
  && \quad {}
  + \iO \Alpha(\phi(t))
  \leq \bigl( \mu(t) + f(t) ,  \phi(t)-v \bigr)
  + \iO \Alpha(v)
  \non
  \\
  && \quad \hbox{for every $v\in\VB{\sigma_0}$\, and \aat}.
  \label{newsecondaz}
\Esist
{Unfortunately, this \newjuerg{does not seem to be} true, in general,
and we prove the following:
\Beq
  \hbox{\sl the variational inequality \eqref{secondaz} with $v$ varying in $\VB{\sigma_0}$ is equivalent to \eqref{newsecondaz}\/.}
  \label{replsecondaz}
\Eeq
To this aim, it suffices to apply the lemma
with the same $\Gam_i$ as before and obvious $u$ and~$g$,
but with $V=\VB{\sigma_0}$ and $a$ defined by $a(u,v):=(u-Pu,v)-\Lpi(u,v)$ for $u,v\in\VB{\sigma_0}$.}

\step
{Conclusion of the proof}

{In view of \eqref{replseconda} and \eqref{replsecondaz}, 
our aim is first to verify \eqref{newsecondaz}
by starting from \eqref{newseconda} (implied by~\eqref{seconda}),
while \eqref{secondaz}, as it is, will be proved at the end by accounting for~\eqref{ass-pier}}.
However, the \lhs\ of \eqref{newseconda} contains the quadratic term associated to the map $v\mapsto-\Lpi\iO|v|^2$.
This term is unpleasant since {the related} map is concave.
To get {rid} of it, we adapt the procedure introduced in~\cite{CG1} to the present case.
We set, for convenience,
\Beq
  \kappa := \frac \Lpi \tau \,, \quad
  \rhos(t) := e^{-\kappa t} \phis(t)
  \aand
  \rho(t) := e^{-\kappa t} \phi(t)
  \quad \aat ,
  \label{notation}
\Eeq
and we notice that $w\mapsto\smash{\intQ e^{-2\kappa t}w^2}$ is the square of an equivalent norm in~$\LQ2$.
At this point, we pick an arbitrary $v\in\L2{\VB{\sigma_0}}$, write \eqref{newseconda} by taking $v(t)$ as test function,
multiply by $e^{-2\kappa t}$, and integrate over~$(0,T)$.
We obtain
\Bsist
  && \ioT \tau \bigl(
    e^{-\kappa t} \bigl( \dt\phis(t) - \kappa \phis(t) \bigr) ,
    e^{-\kappa t} (\phis-v)(t)
  \bigr) \, dt
  + \ioT e^{-2\kappa t} \norma{B^\sigma \phis(t)}^2 \, dt
  \non
  \\
  && \quad {}
  - \ioT e^{-2\kappa t} \bigl( B^\sigma \phis(t) , B^\sigma v(t) \bigr) \, dt
  + \intQ e^{-2\kappa t} \, \Alpha(\phis)
  \non
  \\
  && \leq \ioT e^{-2\kappa t} \bigl( \mus(t) + \fs(t) , (\phis-v)(t) \bigr) \, dt
  + \intQ e^{-2\kappa t} \, \Alpha(v) \,.
  \label{intseconda}
\Esist
Well, we want to take the limit as $\sigma$ tends to zero in this inequality.
As for the first term on the \lhs , we have that
\Bsist
  && \ioT \tau \bigl(
    e^{-\kappa t} \bigl( \dt\phis(t) - \kappa \phis(t) \bigr) ,
    e^{-\kappa t} (\phis-v)(t)
  \bigr) \, dt
  = \ioT {\tau} \bigl( \dt\rhos(t) , \rhos(t) - e^{-\kappa t} v(t) \bigr) \, dt
  \non  
  \\
  && = \frac \tau 2 \, \norma{\rhos(T)}^2
  - \frac \tau 2 \, \norma\phiz^2
  - \ioT \tau \bigl( \dt\rhos(t) , e^{-\kappa t} v(t) \bigr) \, dt \,.
  \non
\Esist
By observing that $\rhos$ converges to $\rho$ weakly in $\H1H$, 
thus weakly in $\C0H$, so that $\rhos(T)$ converges to $\rho(T)$ weakly in~$H$,
we {therefore} have that
\Bsist
  && \liminfs  \ioT \tau \bigl(
    e^{-\kappa t} \bigl( \dt\phis(t) - \kappa \phis(t) \bigr) ,
    e^{-\kappa t} (\phis-v)(t)
  \bigr) \, dt
  \non
  \\
  && \geq \frac \tau 2 \, \norma{\rho(T)}^2
  - \frac \tau 2 \, \norma\phiz^2
  - \ioT \tau \bigl( \dt\rho(t) , e^{-\kappa t} v(t) \bigr) \, dt
  \non
  \\
  && = \ioT \tau \bigl(
    e^{-\kappa t} \bigl( \dt\phi(t) - \kappa \phi(t) \bigr) ,
    e^{-\kappa t} (\phi-v)(t)
  \bigr) \, dt \,.
  \non
\Esist
Next, by \eqref{basic} and the lower semicontinuity of the norms, we have that
\Beq
  \liminfs \ioT e^{-2\kappa t} \norma{B^\sigma \phis(t)}^2 \, dt
  \geq \ioT e^{-2\kappa t} \norma{(\phi-P\phi)(t)}^2 \, dt \,.
  \non
\Eeq
By also recalling \eqref{convBsv}, we can write
\Beq
 \lims \ioT e^{-2\kappa t} \bigl( B^\sigma \phis(t) , B^\sigma v(t) \bigr) \, dt
  = \ioT e^{-2\kappa t} \bigl( (\phi-P\phi)(t) , (v-Pv)(t) \bigr) \, dt \,.
  \non
\Eeq
By taking the difference, we deduce that
\Bsist
  && \liminfs \Bigl(
    \ioT e^{-2\kappa t} \norma{B^\sigma \phis(t)}^2 \, dt
    - \ioT e^{-2\kappa t} \bigl( B^\sigma \phis(t) , B^\sigma v(t) \bigr) \, dt
  \Bigr)
  \non
  \\
  && \geq \ioT e^{-2\kappa t} \norma{(\phi-P\phi)(t)}^2 \, dt
  - \ioT e^{-2\kappa t} \bigl( (\phi-P\phi)(t) , (v-Pv)(t) \bigr) 
  \non
  \\
  && = \ioT e^{-2\kappa t} \bigl( (\phi-P\phi)(t) , (\phi-P\phi)(t) - (v-Pv)(t) \bigr) \, dt
  \non
  \\
  && = \ioT e^{-2\kappa t} \bigl( (\phi-P\phi)(t) , (\phi-v)(t) \bigr) \, dt \, {,}
  \non
\Esist
the last equality being due to the othogonality between $(\phi-P\phi)(t)\in(\ker B)^\perp$ and $(P\phi-Pv)(t)\in\ker B$.
Moreover, by observing that the functional $w{{}\mapsto{}}\intQ e^{-2\kappa t}\Alpha(w)$ is lower semicontinuous on~$\LQ2$,
and recalling that $\phis$ converges to $\phi$ weakly in~$\LQ2$,
we deduce that
\Beq
  \liminfs \intQ e^{-2\kappa t} \Alpha(\phis)
  \geq \intQ e^{-2\kappa t} \Alpha(\phi) \,.
  \non
\Eeq
This ends the treatment of the terms on the \lhs\ of \eqref{intseconda}.
Concerning the \rhs, we have to overcome the difficulty due to the coupling between $\mus$ and~$\phis$.
To this end, we introduce the notation
\Beq
  (1*w)(t) := \iot w(s) \, ds
  \quad \hbox{for every $w\in\L2H$ and $t\in[0,T]$}
  \non
\Eeq
and deduce from \eqref{sprima} that
\Beq
  \phis + A^{2r} (1*{{}\mus{}}) = \phiz \,.
  \non
\Eeq
Hence, we have that
\Bsist
  && \ioT e^{-2\kappa t} \bigl( \mus(t) , (\phis-v)(t) \bigr) \, dt
  = \ioT e^{-2\kappa t} ( \mus(t),\phiz) \, dt
  \non
  \\
  && \quad {}
  - \ioT e^{-2\kappa t} \bigl( A^r \mus(t) , A^r (1*\mus)(t) \bigr) \, dt
  - \ioT e^{-2\kappa t} \bigl( \mus(t) , v(t) \bigr) \, dt \,.
  \non 
\Esist
Now, from \eqref{convmu} we deduce that $1*\mus$ converges to $1*\mu$ weakly in $\H1{\VA{2r}}$.
{Since the embedding $\H1{\VA{2r}}\subset\L2{\VA r}$ is compact, we infer that
$1*\mus$ converges to $1*\mu$ strongly in $\L2{\VA r}$.}
{In view of \eqref{primaz} and \eqref{cauchyz}, we} deduce~that
\Bsist
  && \lims \ioT e^{-2\kappa t} \bigl( \mus(t) , (\phis-v)(t) \bigr) \, dt
  = \ioT e^{-2\kappa t} ( \mu(t),\phiz) \, dt
  \non
  \\
  && \quad {}
  - \ioT e^{-2\kappa t} \bigl( A^r \mu(t) , A^r (1*\mu)(t) \bigr) \, dt
  - \ioT e^{-2\kappa t} \bigl( \mu(t) , v(t) \bigr) \, dt \,.
  \non 
  \\
  && = \ioT e^{-2\kappa t} \bigl( \mu(t) , (\phi-v)(t) \bigr) \, dt \,.
  \non
\Esist
Finally, by recalling \eqref{hpf}, we see that the term involving $\fs$ and the last one of \eqref{intseconda} do not give any trouble.
Therefore, we conclude that
\Bsist
  && \ioT \tau \bigl(
    e^{-\kappa t} \bigl( \dt\phi(t) - \kappa \phi(t) \bigr) ,
    e^{-\kappa t} (\phi-v)(t)
  \bigr) \, dt
  \non
  \\
  && \quad {}
  + \ioT e^{-2\kappa t} \bigl( (\phi-P\phi)(t) , (\phi-v)(t) \bigr) \, dt
  + \intQ e^{-2\kappa t} \, \Alpha(\phis)
  \non
  \\
  && \leq \ioT e^{-2\kappa t} \bigl( \mu(t) + f(t) , (\phi-v)(t) \bigr) \, dt
  + \intQ e^{-2\kappa t} \, \Alpha(v)\,, 
  \label{intsecondaz}
\Esist
and this holds for every $v\in\L2{\VB{\sigma_0}}$.
{On the other hand, \eqref{intsecondaz} is equivalent~to}
\Bsist
  && \tau \bigl(
    e^{-\kappa t} \bigl( \dt\phi(t) - \kappa \phi(t) \bigr) ,
    e^{-\kappa t} (\phi(t)-v)(t)
  \bigr) 
  \non
  \\
  && \quad {}
  + e^{-2\kappa t} \bigl( (\phi-P\phi)(t) , \phi(t) - v \bigr) 
  + \intQ e^{-2\kappa t} \, \Alpha(\phis)
  \non
  \\
  && \leq e^{-2\kappa t} \bigl( \mu(t) + f(t) , \phi(t) - v \bigr) 
  + \iO e^{-2\kappa t} \, \Alpha(v) 
  \non
\Esist
{\aat\ and every $v\in\VB{\sigma_0}$.
By multiplying by $e^{2\kappa t}$ and recalling that $\kappa=\Lpi/\tau$, 
we obtain \eqref{newsecondaz} as claimed.}
{Recalling \eqref{replsecondaz}, we have proved that the variational inequality \eqref{secondaz}
is satisfied for every test function $v\in\VB{\sigma_0}$.
At this point, we account for \eqref{ass-pier}, not yet used up to now,
and show that \eqref{secondaz} actually holds for every $v\in H$.
To this end, for a given $v\in H$ with $\Beta(v)\in\Luno$ without loss of generality,
it suffices to take a sequence $\{\vn\}$ given by \eqref{ass-pier}, test \eqref{newsecondaz} by~$\vn$ and let $n$ tend to infinity.
One obtains \eqref{newsecondaz} for~$v$ without any trouble.
This completes the proof.}

{%
\Brem
\label{No-ass-pier-bis}
Going back to the above proof, 
one justifies what has been announced in Remark~\ref{No-ass-pier}:
if \eqref{ass-pier} is not assumed, one anyway arrives at the variational inequality \eqref{secondaz}
required for every $v\in\VB{\sigma_0}$ instead of \newjuerg{for} every $v\in H$.
Indeed, \eqref{ass-pier} has been only used at the end,
in order to extend to any $v\in H$
the validity of \eqref{secondaz} already proved for test functions $v\in\VB{\sigma_0}$.
\Erem
}%


\section{The limiting problem}
\label{LIMITING}
\setcounter{equation}{0}

In this section, we prove {Theorem \ref{Partuniqz} and Theorem}~\ref{Eqsecondaz}.
As far as the former is concerned,
some preliminaries are needed. 
We refer to \cite[Sect.~3]{CGS18} for more details.
We set
\Beq
  \VA{-r} := \bigl( \VA r \bigr)^*
  \aand
  \norma\cpto_{A,-r}
  := \hbox{the dual norm of $\norma\cpto_{A,r}$}\,, 
  \non
\Eeq
and we use the symbol $\<\cpto,\cpto>_{A,r}$ for the duality pairing between  $\VA{-r}$ and~$\VA r$.
It is understood that $H$ is identified with a subspace of $\VA{-r}$ in the usual way, i.e., in order that
$\<v,w>_{A,r}=(v,w)$ for every $v\in H$ and $w\in\VA r$.
Moreover, we introduce the subspaces $\Vz{\pm r}$ of $\VA{\pm r}$ by setting
\Bsist
  && \Vz r := \VA r
  \aand
  \Vz{-r} := \VA{-r}
  \quad \hbox{if $\lambda_1>0$},
  \non
  \\
  && \Vz r := \{v\in \VA r :\ \mean v=0\}
  \aand
  \Vz{-r} := \{\psi \in \VA{-r} :\ \<\psi,1>_{A,r}=0 \}
  \quad \hbox{if $\lambda_1=0$} \,.
  \non
\Esist
Next, we define $\Az{2r}:\Vz r\to\VA{-r}$ by the formula
\Beq
  \< \Az{2r} v , w >_{A,r} 
  = (A^r v , A^r w)_{A,r}
  \quad \hbox{for every {$v\in\Vz r$} and $w\in\VA r$}.
  \non
\Eeq
It turns out that the range of $\Az{2r}$ is $\Vz{-r}$
and that $\Az{2r}$ is an isomorphism between $\Vz r$ and~$\Vz{-r}$.
Thus, we can set $\Az{-2r}:=(\Az{2r})^{-1}$ 
and obtain an isomorphism between $\Vz{-r}$ and~$\Vz r$.
It also turns out that
\Beq
  \bigl( A^r \Az{-2r} \psi , A^r v )
  = \< \psi,v >_{A,r}
  \quad \hbox{for every $\psi\in\Vz{-r}$ and $v\in\VA r$}.
  \label{ok}
\Eeq
Finally, the following formula holds true:
\Beq
  \< \dt\psi , \Az{-2r} \psi >_{A,r}
  = \frac 12 \, \frac d {dt} \,\norma\psi_{A,-r}^2
  \quad \hbox{\aet, \ for every $\psi\in\H1{\Vz{-r}}$}.
  \non
\Eeq
In particular,
\Beq
  \iot \< \dt\psi(s) , \Az{-2r} \psi(s) >_{A,r} \, ds
  \geq 0
  \quad \hbox{for every $\psi\in\H1{\Vz{-r}}$ with $\psi(0)=0$}.
  \label{nonneg}
\Eeq

\step
Proof of Theorem \ref{Partuniqz}

{We just prove the continuous dependence part, 
since uniqueness for the first component follows as a consequence.}
We set, for convenience,
$f:=f_1-f_2$, $\phi:=\phi_1-\phi_2$, and $\mu:=\mu_1-\mu_2$.
Now, we write equation \eqref{primaz} at the time $s$ for these solutions and take the difference.
Then, we test {the resulting identity by} $\,v=\Az{-2r}\phi(s)${, where we observe} that $\phi(s)\in\Vz{-r}$,
since $\phi\in\C0H$ {by \eqref{regsoluzz}} and $\mean\phi(s)=0$ if $\lambda_1=0$ by the conservation property~\eqref{conservation},
so that $v$ is a well-defined element of~$\VA r$.
Moreover, {we have that} $\Az{-2r}\phi\in\L\infty{\VA r}${.}
Integrating over  $(0,t)$ with respect to~$s$, where $t\in(0,T)$ is arbitrary,
we obtain the identity
\Beq
  \iot \< \dt \phi(s) , \Az{-2r} \phi(s) >_{A,r} \, ds
  + \iot \bigl( A^r \mu(s) , A^r \Az{-2r}\phi(s) \bigr) \, ds
  = 0 \,.
  \non
\Eeq
Now, the first term on the \lhs\ is nonnegative by~\eqref{nonneg}.
Hence, by also noting that $\mu\in\L2{\VA r}$ and applying~\eqref{ok}, 
we deduce that
\Beq
  \iot (\phi(s) , \mu(s)) \, ds 
  \leq 0 \,.
  \label{testprimaz}
\Eeq
At the same time, we write \eqref{secondaz} for $f_i$ and $(\phi_i,\mu_i)$, $i=1,2$,
test them by $\phi_2$ and~$\phi_1$, respectively, 
add the resulting inequalities to each other, and integrate over $(0,t)$ as before.
Then, the terms involving $\Beta$ cancel out. 
By denoting by $I$ the identity map of $H$ and rearranging,
we have that
\Bsist
  && \frac \tau 2 \, \norma{\phi(t)}^2 
  + \iot \bigl( (I-P)\phi(s) , \phi(s) \bigr) \, ds
  \non
  \\
  && \leq \iot \bigl( f(s) + \mu(s) , \phi(s) \bigr) \, ds\,
  - \iot \bigl( \pi(\phi_1(s)) - \pi(\phi_2(s)) , \phi(s) \bigr) \, ds \,.
  \label{testsecondaz}
\Esist
We observe that $I-P$ is the projection operator on the orthogonal subspace $(\ker B)^\perp$.
It follows that $((I-P)v,v)=((I-P)v,(I-P)v)\geq0$ for every $v\in H$,
so that the second term on the \lhs\ of \eqref{testsecondaz} is nonnegative.
By adding \eqref{testprimaz} and \eqref{testsecondaz} to each other, and accounting for this observation, 
an obvious cancellation, the \Lip\ continuity of~$\pi$ and the Schwarz and Young inequalities,
we deduce that 
\Beq
  \frac \tau 2 \, \norma{\phi(t)}^2 
  \leq \frac 14 \iot \norma{f(s)}^2 \, ds
  + (1+\Lpi) \iot \norma{\phi(s)}^2 \, ds.
  \non
\Eeq 
By applying the Gronwall lemma,
we conclude that the desired estimate \eqref{contdep} holds true
with a constant $C_{cd}$ as in the statement.\QED

\medskip

{Finally, we prove Theorem~\ref{Eqsecondaz}.
The proof we give is based on the study of the auxiliary problem of finding $\bphi\in\H1H$ satisfying       
\Bsist
  && \tau \bigl( \dt\bphi(t) , \bphi(t) - v \bigr)
  + \bigl( \bphi(t) - P \bphi(t) , \bphi(t)-v \bigr)
  \non
  \\
  && \quad {}
  + \iO \Beta(\bphi(t))
  + \bigl( \pi(\bphi(t)) - \pi(0) ,  \bphi(t)-v \bigr)
  \non
  \\
  && \leq \bigl( g(t) ,  \bphi(t)-v \bigr)
  + \iO \Beta(v)
  \quad \hbox{for every $v\in H$ and \aat},
  \qquad
  \label{secondaaux}
  \\
  && \bphi(0) = \bphiz, 
  \label{cauchyaux}
\Esist
\Accorpa\Pblaux secondaaux cauchyaux
for given
\Beq
  g \in \L2H
  \aand
  \bphiz \in H \,.
  \label{hpdatiaux}
\Eeq
We have subtracted the constant $\pi(0)$ to ${\pi}(\bphi(t))$ in \eqref{secondaaux}
in order to use the inequality $|\pi(s)-\pi(0{)}|\leq\Lpi\,|s|$ for $s\in\erre$
without any additive constant.
This is needed in the {sequel}, indeed.
Since $\Beta$ is convex, $P$~is linear and $\pi$ is \Lip\ continuous,
this problem has a unique solution~$\bphi$
provided that the initial datum also satisfies
\Beq
  \Beta(\bphiz) \in \Luno \,.
  \label{hpB}
\Eeq
In the forthcoming Lemma~\ref{Aux}, we prove a regularity result
by applying a particular case of \cite[Sect.~I, Thm.~2]{Tartar}
which we present here in the form of a lemma.}

\Blem
\label{Tartar}
Let $\calA_0$, $\calA_1$, $\calB_0$ and $\calB_1$ be four Banach spaces
with the continuous embeddings $\calA_0\subset\calA_1$ and $\calB_0\subset\calB_1$,
and let $\calT:\calA_1\to\calB_1$ be a nonlinear operator satisfying $\calT v\in\calB_0$ for every $v\in\calA_0$.
Assume that
\Bsist
  && \norma{\calT u - \calT v}_{\calB_1}
  \leq C_1 \, \norma{u-v}_{\calA_1}
  \quad \hbox{for every $u,v\in\calA_1$},
  \label{hpC1}
  \\
  && {\norma{\calT v}_{\calB_0}
  \leq C_2 \norma v_{\calA_0}
  \quad \hbox{for every $v\in\calA_0$}},
  \label{hpC2}
\Esist
for some positive constants $C_1$ and~$C_2$.
Then, for every $\theta\in(0,1)$ and $p\in[1,+\infty]$, we have~that
\Bsist
  && \calT v \in (\calB_0,\calB_1)_{\theta,p}
  \aand
  \norma{\calT v}_{(\calB_0,\calB_1)_{\theta,p}}
  \leq C C_1^\theta C_2^{1-\theta} \, \norma v_{(\calA_0,\calA_1)_{\theta,p}}
  \non
  \\
  && \quad \hbox{for every $v\in(\calA_0,\calA_1)_{\theta,p}$},
  \label{tartar}
\Esist
with a constant $C$ that does not depend on~$\calT$.
\Elem

In the above lemma, the symbol $\norma\cpto_X$ stands for the norm in the generic Banach space~$X$.
The same convention is followed in the rest of the section, where $\norma\cpto_X$ also denotes the norm in the power~$X^3$
(however, we keep the short notation $\norma\cpto$ without indices if $X=H$).
Moreover, $(X,Y)_{\theta,p}$ is the real interpolation space between the Banach spaces $X$ and~$Y$
(for basic definitions and properties see, e.g., \cite[Sect.~1.1]{Lunardi}).

\Blem
\label{Aux}
Let the general assumption on the structure be {fulfilled}  
and assume that the data $g$ and $\bphiz$ satisfy
\Beq
  g \in \L2{\Hx\eta}
  \aand 
  \bphiz \in \Hx\eta 
  \label{hpdatieta}
\Eeq
for some $\eta\in\pier{{}(0,1]}$, as well as \eqref{hpB}.
Then, the solution $\bphi$ to problem \Pblaux\ enjoys the further regularity
\Beq
  \bphi \in \L2{\Hx\eta},
  \label{regbphiaux}
\Eeq
and there exists some $\xi$ satisfying
\Bsist
  && \xi \in \L2H
  \aand
  \xi \in \beta(\bphi)
  \quad \aeQ \,,
  \label{regxiaux}
  \\
  && \tau \dt\bphi + \bphi - P\bphi + \xi + \pi(\bphi) - \pi(0) = g
  \quad \aeQ \,.
  \label{eqsecondaaux}  
\Esist
\Elem

\Bdim
By still denoting by $\Betal$ and $\betal$ the Moreau--Yosida approximations of $\Beta$ and~$\beta$, respectively,
we introduce the approximating problem of finding $\phil\in\H1H$ that satisfies
\Beq
  \tau \dt\phil + \phil - P \phil + \betal(\phil) + \pi(\phil) - \pi(0) = g
  \quad \aeQ 
  \label{secondaauxl}
\Eeq
and the initial condition \eqref{cauchyaux}.
For any data satisfying \eqref{hpdatiaux} (while \eqref{hpB} is not needed here)
also this problem has a unique solution~$\phil$.
We perform some a~priori estimates.
As usual, the symbol $c$ stands for possibly different constants.
In this proof, the values of $c$ can only depend on $\tau$, $\pi$, $\Omega$, $T$
and the eigenfunctions $e'_j$ associtated to the zero eigenvalues of~$B$ (if~any).
In particular, they do not depend on~$\lambda$, nor on the data of problem \Pblaux.
Symbols like $C$ and $C_i$ denote particular values of $c$ we want to refer~to.
The first three estimates we perform are in the direction of the inequalities \eqref{hpC1} and \eqref{hpC2}
which we want to satisfy with a suitable choice of the spaces and the operator.
For this reason, they are obtained under different regularity assumptions on the data.

\step
First a priori estimate

Let $g_i$ and $\bphi_{0,i}$, $i=1,2$, be two choices of the data satisfying \eqref{hpdatiaux}
and let $\bphi_{\lambda,i}$ be the corresponding solutions to the approximating problem.
We set for brevity 
$\phil:=\bphi_{\lambda,1}-\bphi_{\lambda,2}$, $g:=g_1-g_2$ and $\bphiz:=\bphi_{0,1}-\bphi_{0,2}$.
We write \eqref{secondaauxl} for both solutions,
take the difference and multiply it by~$\phil$.
Then, we integrate over~$Q_t$.
We obtain~that
\Bsist
  && \frac \tau 2 \iO |\phil(t)|^2
  + \intQt |\phil|^2
  + \intQt \bigl( \betal(\bphi_{\lambda,1}) - \betal(\bphi_{\lambda,2}) \bigr) \phil
  \non
  \\
  && = \frac \tau 2 \iO |\bphiz|^2
  + \intQt g \, \phil
  + \intQt (P\phil) \phil
  - \intQt \bigl( \pi(\bphi_{\lambda,1}) - \pi(\bphi_{\lambda,2}) \bigr) \phil
  \non
\Esist
Since $\betal$ is monotone,
all of the terms on the \lhs\ are nonnegative.
By estimating the \rhs\
on account of the \Lip\ continuity of $\pi$ and the Schwarz and Young inequalities,
and then applying the Gronwall lemma, we easily conclude that
\Beq
  \norma{\bphi_{\lambda,1}-\bphi_{\lambda,2}}_{\L\infty H}
  \leq C_{1,\infty} \, \bigl( \norma{g_1-g_2}_{\L2H} + \norma{\bphi_{0,1}-\bphi_{0,2}} \bigr) \,.
  \label{stima1infty}
\Eeq
It trivially follows that
\Beq
  \norma{\bphi_{\lambda,1}-\bphi_{\lambda,2}}_{\L2 H}
  \leq C_1 \, \bigl( \norma{g_1-g_2}_{\L2H} + \norma{\bphi_{0,1}-\bphi_{0,2}} \bigr) \,.
  \label{stima1}
\Eeq

\step 
Second a priori estimate

We assume \eqref{hpdatiaux} on the data.
By multiplying \eqref{secondaauxl} by $\phil$ and integrating over~$Q_t$,
we obtain~that
\Bsist
  && \frac \tau 2 \iO |\phil(t)|^2
  + \intQt |\phil|^2
  + \intQt \betal(\phil) \phil
  \non
  \\
  && = \frac \tau 2 \iO |\bphiz|^2
  + \intQt g \phil
  + \intQt (P\phil) \phil
  + \intQt \bigl( \pi(\phil) - \pi(0) \bigr) \phil \,.
  \non
\Esist
All of the terms on the \lhs\ are nonnegative since 
$\betal$~is monotone and $\betal(0)=0$.
{If we estimate} the \rhs\ by using the \Lip\ continuity of $\pi$ and the Schwarz and Young inequality,
we immediately deduce~that
\Beq
  \norma\phil_{\L\infty H} \leq c \, ( \norma g_{\L2H} + \norma\bphiz ) \,.
  \label{stimaHinfty}
\Eeq

\step
Third a priori estimate

We set $V:=\Huno$ for brevity and assume that 
the data satisfy $g\in\L2V$ and $\bphiz\in V$.
Before going on, we make an observation.
Assume first that $\ker B=\{0\}$.
Then $P=0$ and \eqref{secondaaux} is an ordinary differential equation
where the space variable is just a parameter.
{In} the opposite case, the presence of the nonlocal operator $P$ could be unpleasant.
However, we are reduced to the same situations as before by moving the term $P\phil$ to the \rhs\ and treating it as a datum.
More precisely, in this case, $\ker B$~has a finite dimension $m>0$ and \newjuerg{is} spanned
by the first $m$ eigenfunctions (those corresponding to the zero eigenvalues).
{Since every eigenfunction of $B$ belongs to the domain $\VB n$ of~$B^n$ for every $n\in\enne$
and we are assuming~\eqref{newhpsuppl}, the eigenfunctions (we are interested in) belong to $V$,
and we have the identities}
\Beq
  Pv = \somma j1m (v,e'_j) e'_j
  {\aand}
  \nabla Pv = \somma j1m (v,e'_j) \nabla e'_j 
  \quad \hbox{for every $v\in H$} \,.
  \label{reprP}
\Eeq
Namely, we have that $Pv\in V$ even though $v$ only belongs to~$H$.
Therefore, in any case, the solution $\phil$ enjoys some space regularity.
Precisely, it belongs to $\L2V$ as well as its time derivative and we have that
\Beq
  \tau \, \dt\nabla\phil
  + \nabla\phil
  + \betal'(\phil) \nabla\phil
  + \pi'(\phil) \nabla\phil
  = \nabla g
  + \nabla P\phil 
  \quad \aeQ \,.
  \non
\Eeq
By multiplying this equation by $\nabla\phil$ and integrating over~$Q_t$,
we obtain~that
\Bsist
  && \frac \tau 2 \iO |\nabla\phil(t)|^2
  + \intQt |\nabla\phil|^2
  + \intQt \betal'(\phil) |\nabla\phil|^2
  \non
  \\
  && = \frac \tau 2 \iO |\nabla\bphiz|^2
  + \intQt \nabla g \cdot \nabla\phil
  + \intQt (\nabla P\phil) \cdot \nabla\phil
  - \intQt \pi'(\phil) |\nabla\phil|^2 \,.
  \non
\Esist
All of the terms on the \lhs\ are nonnegative. The volume integrals 
on the \rhs\newjuerg{, except the one involving~$P$,}
can be easily treated {thanks to the} boundedness of $\pi'$ and the Schwarz and Young inequalities.
If $P=0$\newjuerg{, then} we can apply the Gronwall lemma and obtain an estimate of~$\nabla\phil$.
Recalling \eqref{stimaHinfty}, we conclude~that
\Beq
  \norma\phil_{\L\infty V} \leq C_{2,\infty} \, ( \norma g_{\L2V} + \norma\bphiz_V ) \,.
  \label{stimaVinfty}
\Eeq
We claim that the same estimate holds true even though $\ker B$ is nontrivial.
In this case, we recall the representation formula \eqref{reprP} and apply it to~$\phil$.
By also accounting for standard inequalities, we obtain~that
\Bsist
  && \intQt (\nabla P\phil) \cdot \nabla\phil
  = \intQt \somma j1m (\phil,e'_j) \nabla e'_j \cdot \nabla\phil
  \non
  \\
  && = \somma j1m \iot \Bigl( (\phil(s),e'_j) \iO \nabla e'_j \cdot \nabla\phil(s) \Bigr) \, ds
  \non
  \\
  && \leq \somma j1m \iot \norma{\phil(s)} \, \norma{e'_j} \, \norma{\nabla e'_j} \, \norma{\nabla\phil(s)} \, ds
  \leq c \iot \norma{\phil(s)} \, \norma{\nabla\phil(s)} \, ds
  \non
  \\
  && \leq c \, \norma\phil_{\L2H}^2 + c \intQt |\nabla\phil|^2 \,.
  \non
\Esist
So, it suffices to recall \eqref{stimaHinfty} and apply the Gronwall lemma to obtain \eqref{stimaVinfty} also in this case.
Therefore, \eqref{stimaVinfty} is established and it trivially implies that
\Beq
  \norma\phil_{\L2 V} \leq C_2 \, ( \norma g_{\L2V} + \norma\bphiz_V ) \,.
  \label{stima2}
\Eeq

\step
Interpolation

Now, let the data satisfy~\eqref{hpdatieta} {with $\eta\in(0,1)$}.
We choose
\Bsist
  && \calA_0 := \L2V \times V \,, \quad
  \calA_1 := \L2H \times H \,, \quad
  \non
  \\
  && \calB_0 := \L2V
  \aand
  \calB_1 := \L2H 
  \non
\Esist
and apply Lemma~\ref{Tartar} to the operator $\calT:\calA_1\to\calB_1$ 
that associates to the pair $(g,\bphiz)$ the solution $\phil$ to problem \Pblaux.
Then, \eqref{stima1} and \eqref{stima2} yield \eqref{hpC1} and~\eqref{hpC2}, respectively.
Moreover, by setting $\theta:=1-\eta$, we have that
\Beq
  (\calA_0,\calA_1)_{\theta,2}
  = (\L2V,\L2H)_{\theta,2} \times (V,H)_{\theta,2}
  = \L2{\Hx\eta} \times \Hx\eta
  \non
\Eeq
so that $(g,\bphiz)\in(\calA_0,\calA_1)_{\theta,2}$ by \eqref{hpdatieta}.
It follows that
\Bsist
  && \phil \in (\calB_0,\calB_1)_{\theta,2}
  = \L2{\Hx\eta}
  \aand
  \non
  \\
  && \norma\phil_{\L2{\Hx\eta}}
  \leq C C_1^\theta C_2^{1-\theta} \norma{(g,\bphiz)}_{\L2{\Hx\eta} \times \Hx\eta}
  \label{datartar}
\Esist
with a constant $C$ that does not depend on~$\lambda$.
{Notice that \eqref{datartar} with $\eta=1$ (i.e., $\theta=0$) is~ensured by~\eqref{stima2}.}

\step
Fourth a priori estimate

We are close to the conclusion, and we thus assume that
the data $g$ and $\bphiz$ are as in the statement.
By multiplying \eqref{secondaauxl} by $\dt\phil$, integrating over~$Q_t$, and rearranging,
we have~that
\Bsist
  && \tau \intQt |\dt\phil|^2
  + \frac 12 \iO |\phil(t)|^2
  + \iO \Betal(\phil(t))
  \non
  \\
  && {} = \frac 12 \iO |\bphiz|^2
  + \iO \Betal({\bphiz})
  + \intQt \bigl( g + P\phil - \pi(\phil) + \pi(0) \bigr) \dt\phil \,.
  \non
\Esist
Since $\Betal$ is nonnegative and $\Betal({\bphiz})\leq\Beta({\bphiz})$ \aeO, owing to the Schwarz and Young inequalities and the \Lip\ continuity of $\pi$, 
and accounting for \eqref{hpB} and~\eqref{stimaHinfty},
we infer~that
\Beq
  \norma{\dt\phil}_{\L2H}
  \leq c \bigl( \norma g_{\L2H} + \norma\bphiz + \norma{\Beta({\bphiz})}_{\Luno}^{1/2} \bigr) \,.
  \label{stimadtphil}
\Eeq
A comparison in \eqref{secondaauxl} then yields that
\Beq
  \norma{\betal(\phil)}_{\L2H}  
  \leq c \bigl( \norma g_{\L2H} + \norma\bphiz + \norma{\Beta({\bphiz})}_{\Luno}^{1/2} \bigr) \,.
  \label{stimabetal}
\Eeq

\step
Conclusion

At this point, we let $\lambda$ tend to zero based on \accorpa{datartar}{stimabetal},
the compact embedding $\Hx\eta\subset H$ {for $\eta\in(0,1]$},
and the well-known Aubin--Lions lemma (see, e.g., \cite[Thm.~5.1, p.~58]{Lions}).
We deduce that there exists a pair $(\bphi,\xi)$ such that
\Bsist
  && \phil \to \bphi
  \quad \hbox{weakly star in $\H1H\cap\L2{\Hx\eta}$}
  \non
  \\
  && \quad \hbox{and strongly in $\L2H$}\,,
  \label{convphil}
  \non
  \\
  && \betal(\phil) \to \xi
  \quad \hbox{weakly in $\L2H$}\,,
  \label{convxil}
\Esist
possibly only for a subsequence $\lambda_k\searrow0$. 
Then, $\bphi(0)=\bphiz$, and \eqref{eqsecondaaux} is verified.
Moreover, by also applying, e.g., \cite[Lemma~2.3, p.~38]{Barbu}, 
we infer~that $(\bphi,\xi)$ satisfies the {inclusion in}~\eqref{regxiaux}.
On the other hand, all this implies \eqref{secondaaux} since $\Beta$ is convex,
so that $\bphi$ is the solution to problem \Pblaux.
This completes the proof of the lemma.
\Edim

\step
Proof of Theorem~\ref{Eqsecondaz}

We apply Lemma~\ref{Aux} by choosing
\Beq
  g = \mu + f - \pi(0)
  \aand
  \bphiz = \phiz \,.
  \non
\Eeq
Notice that conditions \eqref{hpdatieta} are satisfied due to \accorpa{hpsuppl}{regsoluz}.
We thus obtain the existence of some $\xi$ satisfing \eqref{regxiaux} and~\eqref{eqsecondaaux}.
The latter reads
\Beq
  \tau \dt\bphi + \bphi - P\bphi + \xi + \pi(\bphi) - \pi(0)
  = \mu + f - \pi(0)
  \quad \aeQ \,.
  \non
\Eeq
But $\phi$ satisfies this equation {(see Remark~\ref{From-ass-pier})}
since $(\phi,\mu)$ is a solution to problem \Pblz\ by assumption,
and this implies \eqref{secondaaux} for $\phi$ since $\Beta$ is convex.
On the other hand, we have that $\phi(0)=\phiz=\bphiz$.
Since the solution $\bphi$ to problem \Pblaux\ is unique, we conclude that $\bphi=\phi$.
Therefore, \accorpa{regxi}{eqsecondaz} are proved.
The last sentence of the statement trivially follows.
\QED

\Brem
\label{No-use-ass-pier}
We observe that in Theorem~\ref{Eqsecondaz} we start from a solution $(\phi,\mu)$ to problem \Pblz\
without using sufficient conditions for the existence of such a solution.
In particular, \eqref{ass-pier} is not accounted for.
We also notice that the argument followed in the above proof 
provides the existence of a unique solution $\phi$ to both 
equation \eqref{eqsecondaz} and the variational inequality \eqref{secondaz}
for a given~$\mu$ without the use of~\eqref{ass-pier}.
\Erem

{%
\Brem
\label{Different}
It is possible to slightly modify the proof of Lemma~\ref{Aux} 
in the application of Lemma~\ref{Tartar} and to obtain different regularity results in Theorem~\ref{Eqsecondaz}.
One can play with the index~$p$ in the interpolation argument, indeed.
If we want to maximize the time regularity, we change the choice of the spaces $\calB_i$
by taking
\Beq
  \calB_0 := \L\infty V
  \aand
  \calB_1 := \L\infty H 
\Eeq
and start from \eqref{stima1infty} and \eqref{stimaVinfty}
in place of \eqref{stima1} and~\eqref{stima2}.
Then, we apply Lemma~\ref{Tartar} still with $\theta=1-\eta$, but with $p=\infty$.
Instead of \eqref{datartar}, we obtain~that
\Bsist
  && \phil \in (\L\infty V,\L\infty H)_{\theta,\infty}
  \aand
  \non
  \\
  && \norma\phil_{(\L\infty V,\L\infty H)_{\theta,\infty}}
  \leq C C_1^\theta C_2^{1-\theta} \norma{(g,\bphiz)}_{\L2{\Hx\eta} \times \Hx\eta}\,,
  \non
\Esist
still with a constant $C$ that does not depend on~$\lambda$.
Then everything can proceed as before.
At the end of the proof of Theorem~\ref{Eqsecondaz}, we arrive at the regularity
\Beq
  \phi \in (\L\infty V,\L\infty H)_{\theta,\infty}
  \label{besov}
\Eeq
for the first component $\phi$ of the solution $(\phi,\mu)$ to problem \Pblz.
We avoid the troubles that may arise with the exponent $\infty$
and do not offer a different representation of the space appearing in~\eqref{besov}.
We just remark that the regularity \eqref{besov} is neither better nor worse than~\eqref{regxi}, 
since it yields some better time regularity \newjuerg{at the expense of} a lower space regularity.
One can prove that $(\L\infty V,\L\infty H)_{\theta,\infty}\subset\L\infty{\Hx{\eta-\eps}}$ for every $\eps>0$
(in~particular, the Aubin--Lions lemma can be applied also in the modified proof of Lemma~\ref{Aux})
so~that the Sobolev type regularity for $\phi$ we can obtain~is
\Beq
  \phi \in \L\infty{\Hx{\eta-\eps}}
  \quad \hbox{for every $\eps>0$} \,.
  \non
\Eeq
\Erem
}%

\Brem
\label{Uniquemuz}
Concerning uniqueness for the second component $\mu$ of the solution to problem \Pblz,
{we can give sufficient conditions in a different direction}.
The situation is similar to the one encountered for problem \Pbl\ and mentioned in Remark~\ref{Uniquemu}.
Let us give some detail.
Assume that $(\phi,\mu_i)$, $i=1,2$, are solutions corresponding to some data $\phiz$ and~$f$
(with the same first component, due to Theorem~\ref{Partuniqz}).
By writing \eqref{primaz} for both solutions and taking the difference,
we immediately obtain that $(A^r(\mu_1-\mu_2),v)=0$ for every $v\in\VA r$ and \aet, that~is
\Beq
  A^r(\mu_1-\mu_2) = 0 \,.
  \label{foruniquemu}
\Eeq
This implies that $\mu_1=\mu_2$ if $\lambda_1>0$.
In the opposite case $\lambda_1=0$, 
we can arrive at the same conclusion under additional conditions,
as we show at once by following the ideas of \cite[Rem.~4.1]{CGS18}.
However, in the present case, the condition we assume on the solutions is difficult to verify, unfortunately.
Suppose that $D(\beta)$ is an open interval, the restriction of $\Beta$ to $D(\beta)$ is a $C^1$ function,
and all {of the values attained by} $\phi$ belong to a compact subinterval $[a,b]\subset D(\beta)$.
Now, choose $\delta_0$ such that the interval $[a-\delta_0,b+\delta_0]$ is contained in~$D(\beta)$.
Then, for an arbitrary $\delta\in(0,\delta_0)$ and \aat,
we can choose $v=\phi(t)-\delta$ (whence $\phi(t)-v=\delta$) and $v=\phi(t)+\delta$ (whence $\phi(t)-v=-\delta$)
in~the variational inequality \eqref{secondaz} written for $(\phi,\mu_1)$ and $(\phi,\mu_2)$, respectively.
Then, by adding the resulting inequalities, we deduce that 
\Beq
  2\iO \Beta (\phi)
  \leq \delta (\mu_1 - \mu_2, \1)
  + \iO \Beta (\phi-\delta)
  + \iO \Beta (\phi+\delta)
  \quad \aet.
  \non 
\Eeq
Division by $\delta$ then yields that
\Beq
  \iO \frac{\Beta (\phi)-\Beta (\phi-\delta)} \delta
  + \iO \frac{\Beta(\phi)-\Beta(\phi+\delta)} \delta
  \leq (\mu_1 - \mu_2, \1) . 
  \non
\Eeq
Taking the limit as $\delta \searrow 0$, we conclude from the Lebesgue dominated convergence theorem that 
\Beq
  0 = \iO \beta(\phi) - \iO \beta(\phi)
  \leq (\mu_1 - \mu_2, \1). 
  \non
\Eeq
Interchanging the roles of $\mu_1$  and $\mu_2$, 
we then infer that $\mean \mu_1 = \mean \mu_2$ \aet. 
By combining this with \eqref{foruniquemu}, we conclude that $\mu_1 = \mu_2$. 
\Erem


\section*{Acknowledgments}
This research was supported by the Italian Ministry of Education, 
University and Research~(MIUR): Dipartimenti di Eccellenza Program (2018--2022) 
-- Dept.~of Mathematics ``F.~Casorati'', University of Pavia. 
In addition, {PC gratefully mentions} some other support 
from the GNAMPA (Gruppo Nazionale per l'Analisi Matematica, 
la Probabilit\`a e le loro Applicazioni) of INdAM (Isti\-tuto 
Nazionale di Alta Matematica), and PC and CG point out their affiliations 
as Research Associates to the IMATI -- C.N.R. Pavia, Italy.


\vspace{3truemm}

\Begin{thebibliography}{10}

\bibitem{AM}
M. Ainsworth, Z. Mao, Analysis and approximation of a fractional
Cahn--Hilliard equation. SIAM J. Numer. Anal. {\bf 55} (2017),
1689-1718.

\bibitem{AkSS1}
G. Akagi, G. Schimperna, A. Segatti,
Fractional Cahn--Hilliard, Allen--Cahn and porous medium equations,
{\em J. Differential Equations} {\bf 261} (2016), 2935-2985.

\bibitem{AkSS2}
G. Akagi, G. Schimperna, A. Segatti,
Convergence of solutions for the fractional Cahn--Hilliard system.
{\em J. Funct. Anal.}  {\bf 276}  (2019), 2663-2715.

\bibitem{Barbu}
V. Barbu,
``Nonlinear Differential Equations of Monotone Type in Banach Spaces'',
Springer, London, New York, 2010.

\pier{\bibitem{BCGG} 
V. Barbu, P. Colli, G. Gilardi, M. Grasselli,
Existence, uniqueness, and longtime behavior
for a nonlinear Volterra integrodifferential equation,
{\em Differential Integral Equations} {\bf 13} (2000), 1233-1262.}

\bibitem{BCST}
E. Bonetti, P. Colli, L. Scarpa, G. Tomassetti,  
Bounded solutions and their asymptotics for a doubly nonlinear Cahn--Hilliard system.
{\em Calc. Var. Partial Differential Equations}, doi:10.1007/s00526-020-1715-9.

\bibitem{Brezis}
H. Brezis,
``Op\'erateurs maximaux monotones et semi-groupes de contractions
dans les espaces de Hilbert'',
North-Holland Math. Stud.
{\bf 5},
North-Holland,
Amsterdam,
1973.

\bibitem{BMH}	
L. Bu, L. Mei, Y. Hou, 
Stable second-order schemes for the space-fractional Cahn--Hilliard and Allen--Cahn equations,
{\em Comput. Math. Appl.}  {\bf 78}  (2019),  3485-3500.

\bibitem{CRW}
C. Cavaterra, E. Rocca, H. Wu,
Long-time dynamics and optimal control of a 
diffuse interface model for tumor growth,
{\em Appl. Math. Optim.}, doi:10.1007/s00245-019-09562-5.

\bibitem{ChFaPr}
R. Chill, E. Fa\v sangov\'a, J. Pr\"uss, 
Convergence to steady state of solutions of the Cahn--Hilliard and Caginalp equations with dynamic boundary conditions, 
{\em Math. Nachr.} {\bf 279} (2006), 1448-1462.

\bibitem{ClGu}
B. Climent-Ezquerra, F.  Guill\'en-Gonz\'alez,  
Convergence to equilibrium of global weak solutions for a Cahn-Hilliard-Navier-Stokes vesicle model,
{\em Z. Angew. Math. Phys.} {\bf  70}  (2019),  Paper No. 125, 27 pp.

\bibitem{CF}
P. Colli, T. Fukao, 
Nonlinear diffusion equations as asymptotic 
limits of Cahn--Hilliard systems, 
{\em J. Differential Equations} {\bf 260} (2016), 6930-6959.

\bibitem{CG1}
P. Colli, G. Gilardi, 
Well-posedness, regularity and asymptotic analyses 
for a fractional phase field system, 
{\em Asympt. Anal.} {\bf 114} (2019), 93-128.

\bibitem{CGS17}
P. Colli, G. Gilardi, J. Sprekels,
On the longtime behavior of a viscous Cahn--Hilliard system 
with convection and dynamic boundary conditions,
{\em J. Elliptic Parabol. Equ.} {\bf 4} (2018), 327-347.

\bibitem{CGS18}
P. Colli, G. Gilardi, J. Sprekels,
Well-posedness and regularity for a generalized fractional Cahn--Hilliard system,
{\em Atti Accad. Naz. Lincei Rend. Lincei Mat. Appl.} {\bf 30} (2019), 437-478.

\bibitem{CGS21bis}
P. Colli, G. Gilardi, J. Sprekels,
Recent results on well-posedness and optimal control for a class of generalized
fractional Cahn--Hilliard systems, {\em Control Cybernet.} {\bf 48} (2019), 153-197.

\bibitem{CGS19}
P. Colli, G. Gilardi, J. Sprekels,
Optimal distributed control of a generalized fractional Cahn--Hilliard system,
{\em Appl. Math. Optim.}, doi:10.1007/s00245-018-9540-7.

\bibitem{CGS24}
P. Colli, G. Gilardi, J. Sprekels, 
Asymptotic analysis of a tumor growth model 
with fractional operators, 
{\em Asymptot. Anal.}, doi:10.3233/ASY-191578.

\bibitem{CGS22}
P. Colli, G. Gilardi, J. Sprekels, 
Longtime behavior for a generalized 
Cahn--Hilliard system with fractional operators,
{{\em Atti Accad. Peloritana Pericolanti Cl. Sci. Fis. Mat. Natur.},
to appear (see also the preprint arXiv:1904.00931~[math.AP] (2019), pp.~1-18).}

\bibitem{Gal3}
C. G. Gal,
Well-posedness and long time behavior of the non-isothermal 
viscous Cahn--Hilliard equation with dynamic boundary conditions,
{\em Dyn. Partial Differ. Equ.} {\bf 5} (2008), 39-67.

\bibitem{Gal2}
C. G. Gal,
Non-local Cahn--Hilliard equations with fractional dynamic boundary,
{\em European J. Appl. Math.} {\bf 28} (2017), 736-788.

\bibitem{GiMiSchi} 
G. Gilardi, A. Miranville, G. Schimperna,
On the Cahn--Hilliard equation with irregular potentials and dynamic boundary conditions,
{\em Commun. Pure Appl. Anal.} {\bf 8} (2009), 881-912.

\bibitem{GiMiSchi2} 
G. Gilardi, A. Miranville, G. Schimperna,
Long-time behavior of the Cahn-Hilliard equation with irregular potentials and dynamic boundary conditions
{\em Chinese Annals of Mathematics} Series B {\bf 31} (2010) 679-712. 

\bibitem{GS}
G. Gilardi, J. Sprekels,
Asymptotic limits and optimal control for the Cahn--Hilliard 
system with convection and dynamic boundary conditions,
{\em Nonlinear Anal.} {\bf 178} (2019), 1-31. 

\bibitem{GrPeSch}
M. Grasselli, H. Petzeltov\'a, G. Schimperna,  
Asymptotic  behavior  of  a  nonisothermal  viscous  Cahn--Hilliard  equation  with  inertial  term,
{\em J. Differential Equations} {\bf 239} (2007), 38-60.

\bibitem{Lions}
J.-L. Lions, ``Quelques M\'ethodes de R\'esolution 
des Probl\`emes aux Limites non Lin\'eaires'', 
Dunod Gauthier-Villars, Paris, 1969. 

\bibitem{Lunardi}
A. Lunardi,
``Interpolation Theory'', third edition,
Lecture Notes Scuola Normale Superiore di Pisa, Series Appunti {\bf 16}, Edizioni della Normale, Pisa, 2018.

\bibitem{MiZe}
A. Miranville, S. Zelik, 
Robust exponential attractors for Cahn--Hilliard type equations with singular potentials, 
{\em Math. Methods Appl. Sci.} {\bf 27} (2004), 545-582.


\bibitem{Tartar}
L. Tartar,
Interpolation nonlin\'eaire et r\'egularit\'e,
{\em J. Funct. Anal.}  {\bf 9}  (1972), 469-489.

\bibitem{WCW}
F. Wang, H. Chen, H. Wang, 
Finite element simulation and efficient algorithm for fractional Cahn--Hilliard equation.
{\em J. Comput. Appl. Math.}  {\bf 356}  (2019), 248-266.

\bibitem{WaWu}
X.-M. Wang, H. Wu, 
Long-time behavior for the Hele--Shaw--Cahn--Hilliard system,
{\em Asymptot. Anal.} {\bf 78} (2012), 217-245.

\bibitem{WuZh}
H. Wu, S. Zheng, 
Convergence to equilibrium for the Cahn--Hilliard equation with dynamic boundary conditions, 
{\em J. Differential Equations} {\bf 204} (2004), 511-531.

\bibitem{YLZ}
H. Ye, Q.  Liu, M. Zhou,
An $L^\infty$ bound for solutions of a fractional Cahn--Hilliard equation,
{\em Comput. Math. Appl.}  {\bf 79}  (2020),  3353-3365.

\bibitem{ZhaoLiu}
X. Zhao, C. Liu,
On the existence of global attractor for 3D viscous Cahn--Hilliard equation,
{\em Acta Appl. Math.} {\bf 138} (2015), 199-212.

\bibitem{Zheng}
S. Zheng,
Asymptotic behavior of solution to the Cahn--Hilliard equation,
{\em Appl. Anal.} {\bf 23} (1986), 165-184.

\End{thebibliography}

\End{document}
